\documentclass[review]{elsarticle}

\usepackage{lineno,hyperref}
\modulolinenumbers[5]

\usepackage{xcolor}
\usepackage{amsmath}
\usepackage{amsfonts}
\usepackage{stmaryrd}
\usepackage{array}
\usepackage{caption}
\newcolumntype{M}[1]{>{\centering\arraybackslash}m{#1}}

\usepackage[ruled,vlined]{algorithm2e}
\journal{Journal of Computational Physics}

\begin{document}

\begin{frontmatter}

\title{Refinement of polygonal grids\\
using Convolutional Neural Networks\\
with applications to polygonal Discontinuous\\
Galerkin and Virtual Element methods}%\tnoteref{mytitlenote}}

% %% Group authors per affiliation:
% \author{P. F. Antonietti\fnref{myfootnote}}
% \ead{paola.antonietti@polimi.it}
% \author{E. Manuzzi\fnref{myfootnote}\corref{mycorrespondingauthor}}
% \ead{enrico.manuzzi@polimi.it}
% \address{Politecnico di Milano, p.zza Leonardo da Vinci, 32, 20133 Milano, Italy}
% \fntext[myfootnote]{P. F. Antonietti and E. Manuzzi are members of INDAM-GNCS.}

% \cortext[mycorrespondingauthor]{Corresponding author}


\author[mymainaddress]{P. F. Antonietti\fnref{myfootnote}}
\ead{paola.antonietti@polimi.it}

\author[mymainaddress]{E. Manuzzi\fnref{myfootnote}\corref{mycorrespondingauthor}}
\ead{enrico.manuzzi@polimi.it}

\address[mymainaddress]{Politecnico di Milano, p.zza Leonardo da Vinci, 32, 20133 Milano, Italy}
\fntext[myfootnote]{P. F. Antonietti and E. Manuzzi are members of INDAM-GNCS.}

\cortext[mycorrespondingauthor]{Corresponding author}

% \tnoteref{mytitlenote}}
\tnotetext[mytitlenote]{Abbreviations: Convolutional Neural Networks (CNNs), Polygonal Discontinuous Galerkin (PolyDG), Virtual  Element  Methods  (VEMs), Mid-Point (MP), CNN-enhanced Mid-Point (CNN-MP) refinement strategy, CNN-enhanced  Reference  Polygon  (CNN-RP) refinement strategy.}

\begin{abstract}
We propose new strategies to handle polygonal grids refinement based on Convolutional Neural Networks (CNNs). We show that CNNs can be successfully employed to identify correctly the ``shape" of a polygonal element so as to design suitable refinement criteria to be possibly employed within adaptive refinement strategies. We propose two refinement strategies that exploit the use of CNNs to classify elements' shape, at a low computational cost. We test the proposed idea considering two families of finite element methods that support arbitrarily shaped polygonal elements, namely Polygonal Discontinuous Galerkin (PolyDG) methods and Virtual Element Methods (VEMs). We demonstrate that the proposed algorithms can greatly improve the performance of the discretization schemes both in terms of accuracy and quality of the underlying grids. Moreover, since the training phase is performed off-line and is independent of the differential model the overall computational costs are kept low.
\end{abstract}


\end{frontmatter}

%\linenumbers

%%%%
\section{Introduction}
{\color{black}Many applications in the fields of engineering and applied sciences, such as fluid-structure interaction problems, flow in fractured porous media, crack and wave propagation problems, are characterized by a strong complexity of the physical domain, possibly involving moving geometries, heterogeneous media, immersed interfaces (such as e.g. fractures) and complex topographies. Whenever classical Finite Element methods are employed to discretize the underlying differential model, the process of grid generation can be the bottleneck of the whole simulation, as computational meshes can be composed only of tetrahedral, hexahedral, or prismatic elements. To overcome this limitation, in the last years there has been a great interest in developing finite element methods that can employ general polygons and polyhedra as grid elements for the numerical discretizations of partial differential equations.}
%In the last years, there has been a great interest in developing polygonal finite element methods for the numerical discretizations of partial differential equations.
We mention the mimetic finite difference method \cite{hyman1997numerical,brezzi2005family,brezzi2005convergence,da2014mimetic}, the hybridizable discontinuous Galerkin methods \cite{cockburn2008superconvergent,cockburn2009superconvergent,cockburn2009unified,cockburn2010projection}, the Polygonal Discontinuous Galerkin (PolyDG) method \cite{hesthaven2007nodal,bassi2012flexibility,antonietti2013hp,cangiani2014hp,antonietti2016review,cangiani2017hp}, the Virtual Element Method (VEM) \cite{beirao2013basic,beirao2014hitchhiker,beirao2016virtual,da2016mixed} and the hybrid high-order method \cite{di2014arbitrary,di2015hybrid,di2015hybrid2,di2016review,di2019hybrid}.
This calls for the need to develop effective algorithms to handle polygonal and polyhedral grids and to assess their quality (see e.g. \cite{attene2019benchmark}).
{\color{black}For a comprehensive overview we refer to the monographs and special issues \cite{da2014mimetic,cangiani2017hp,di2016review,di2021polyhedral,beirao2014hitchhiker,burman2021convergence}.}
Among the open problems, there is the issue of handling polytopic mesh refinement \cite{lai2016recursive,hoshina2018simple,berrone2021refinement}, i.e. partitioning mesh elements into smaller elements to produce a finer grid, and agglomeration strategies \cite{chan1998agglomeration,antonietti2020agglomeration,bassi2012flexibility}, i.e. merging mesh elements to obtain coarser grids. {\color{black}Indeed, during either refinement or agglomeration it is important to preserve the quality of the underlying mesh, because this might affect the overall performance of the method in terms of stability and accuracy. Using a suitable mesh may allow to achive the same accuracy with a much smaller number of degrees of freedom when solving the numerical problem, hence saving memory and computational power. However, since in such a general framework mesh elements may have any shape, there are not well established strategies to achieve effective refinement or agglomeration with a fast, robust and simple approach, contrary to classical Finite Elements.} Moreover, grid agglomeration is a topic quite unexplored, because it is not possible to develop such kind of strategies within the framework of classical Finite Elements.\\
In this work, we propose a new strategy to handle polygonal grid refinement based on Convolutional Neural Networks (CNNs). {\color{black} CNNs are powerful function approximators used in machine learning}, that are particularly well suited for image classification when clearly	defined rules cannot be deduced. Indeed, they have been successfully applied in many areas, especially computer vision \cite{lecun2015deep}. In recent years there has been a great development of machine learning algorithms to enhance and accelerate numerical methods for scientific computing. Examples include, but are not limited to,  \cite{raissi2019physics,raissi2018hidden,regazzoni2019machine,regazzoni2020machine,hesthaven2018non,ray2018artificial}.\\
In this work we show that CNNs can be successfully employed to identify correctly the “shape” of a polygonal element without resorting to any geometric property. This information can then be exploited to apply tailored refinement strategies for different families of polygons. This approach has several advantages:
\begin{itemize}
    \item It helps preserving the mesh quality, since it can be easily tailored for different types of elements.
    \item It can be combined with suitable (user-defined) refinement strategies.
    \item It is independent of the numerical method employed to discretize the underlying differential model.
    \item The overall computational costs are kept low, since the training phase of a CNN is performed off line and it is independent of the differential model at hand.
\end{itemize}
{\color{black}The proposed approach is general and can be extended in three dimension, provided that suitable refinement strategies are available for polyhedra.}
In this paper, we show that CNNs can be used effectively to boost either existing refinement criteria, such as the Mid-Point (MP) strategy, that consists in connecting the edges midpoints of the polygon with its centroid, and we also propose a refinement algorithm that employs pre-defined refinement rules on regular polygons. We refer to these paradigms as CNN-enhanced refinement strategies. To demonstrate the capabilities of the proposed approach we consider a second-order model problem discretized by either PolyDG methods and VEMs and we test the two CNN-enhanced refinement strategies based on polygons' shape recognition. For both the CNNs-enhanced refinement strategies we demonstrate their effectiveness through an analysis of quality metrics and accuracy of the discretization methods.\\ 
\newline
The paper is organized as follows. In Section 2 we show how to classify polygons' shape using a CNN. In Section 3 we present new CNN-enhanced refinement strategies and different metrics to measure the quality of the proposed refinement strategies. In Section 4 we train a CNN for polygons classification. In Section 5 we validate the new refinement strategies on a set of polygonal meshes, whereas in Section 6 we test them when employed with PolyDG and Virtual Element discretizations of a second-order elliptic problem. Finally in Section 7 we draw some conclusions.

\section{Polygon classification using CNNs}
\begin{figure}%[!ht]
    %\centering
    \includegraphics[width=0.24\linewidth]{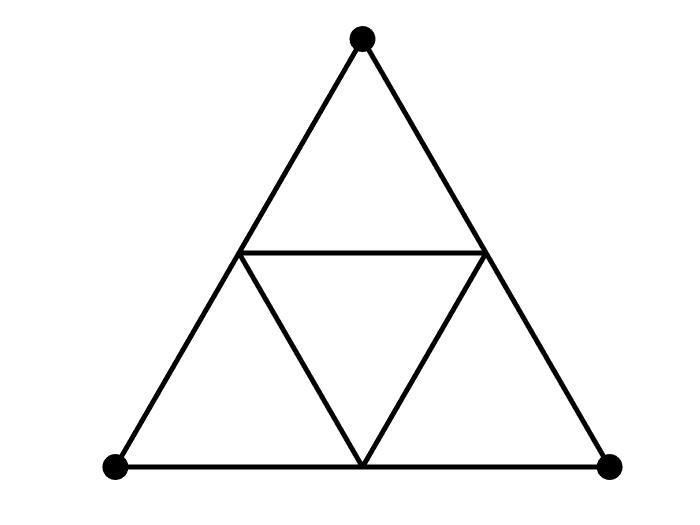}
    \includegraphics[width=0.24\linewidth]{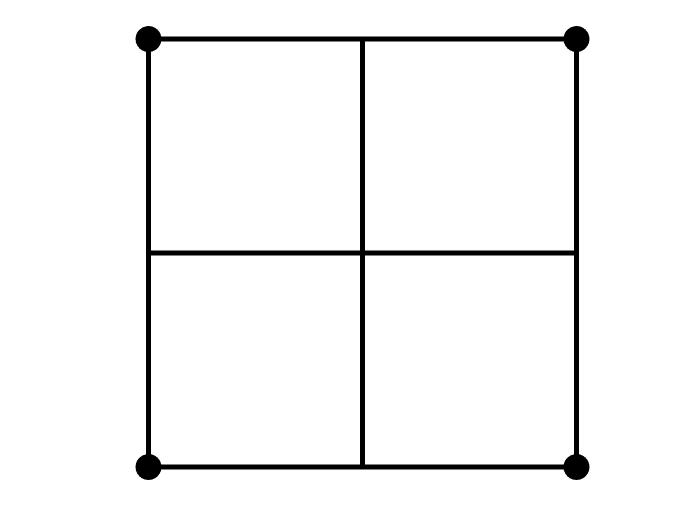}
    \includegraphics[width=0.24\linewidth]{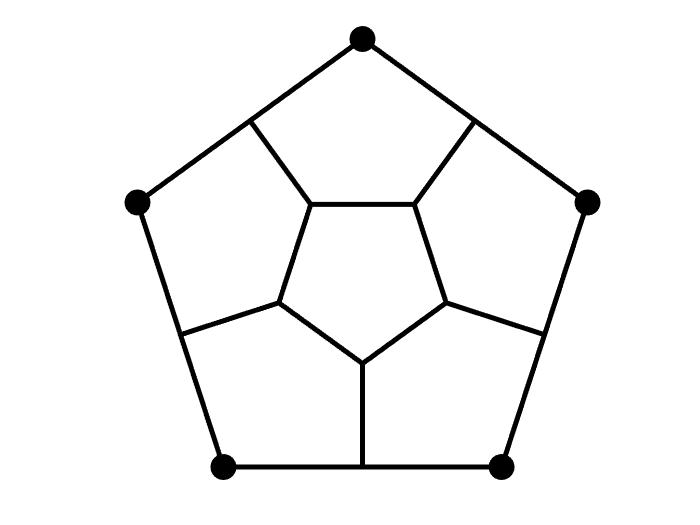}
    \includegraphics[width=0.24\linewidth]{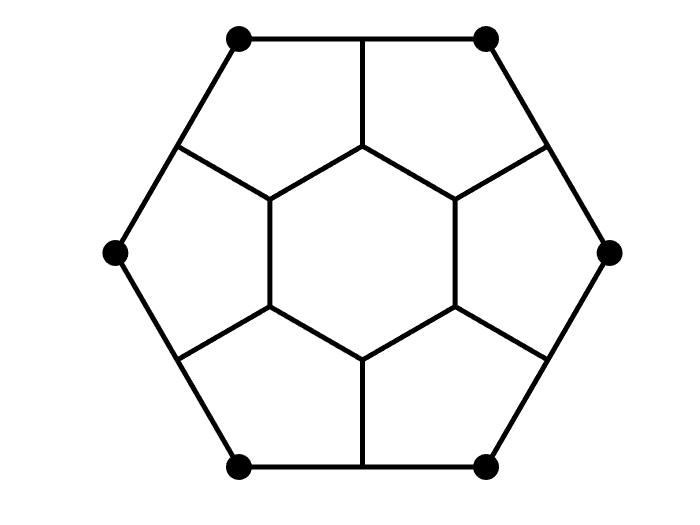}

    \caption[Refinement strategies for regular polygons.]{Refinement strategies for triangular, quadrilateral, pentagonal and hexagonal polygons. The vertices of the original polygon $K$ are marked with black dots.}
    \label{fig:regualar}
\end{figure}

\begin{figure}%[!ht]
    \centering
    \includegraphics[width=0.32\linewidth]{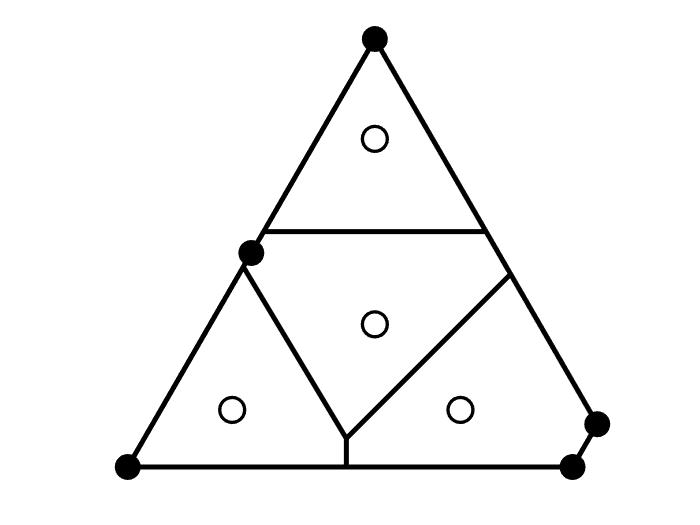}
    \includegraphics[width=0.32\linewidth]{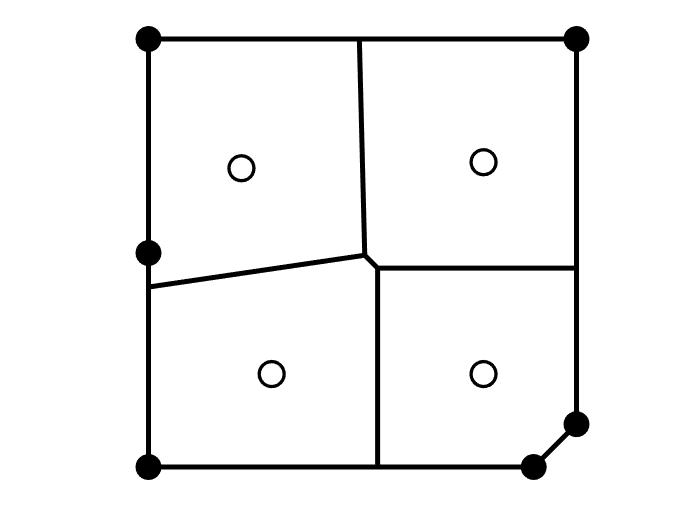}
    \caption[Refinement using Voronoi tessellation.]{Refinement using Voronoi tessellation. The vertices of the original polygon $K$ are marked with full dots, while seeds are marked with empty dots.}
    \label{fig:voronoi ref}
\end{figure}

\begin{figure}[!ht]
    \centering
        \includegraphics[width=0.3\linewidth]{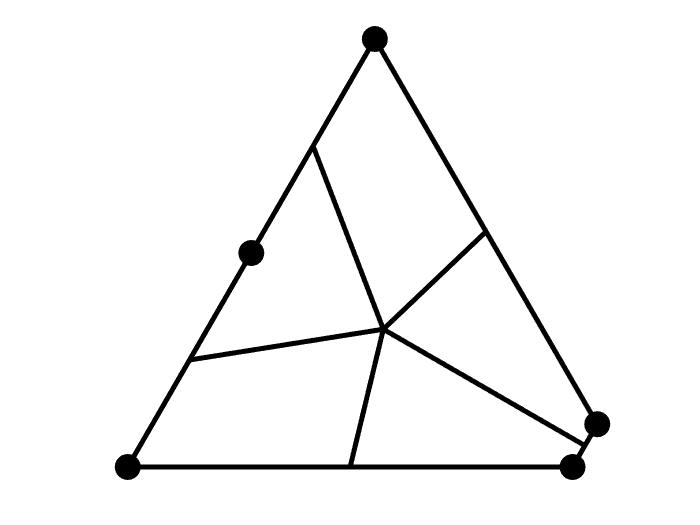}
    \includegraphics[width=0.3\linewidth]{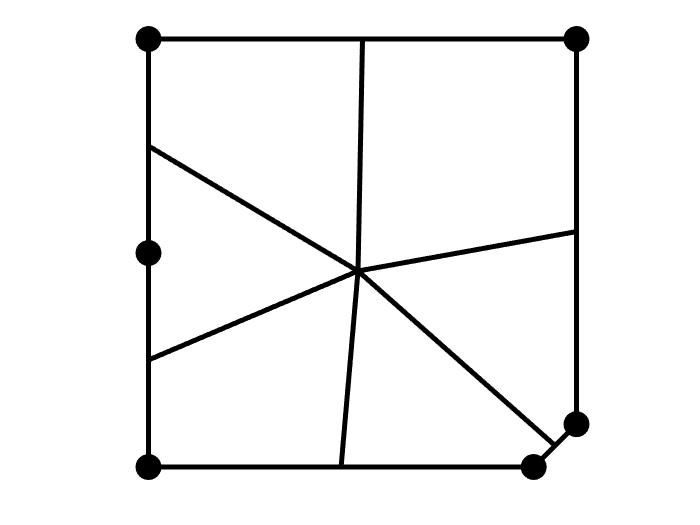}
    
    \includegraphics[width=0.3\linewidth]{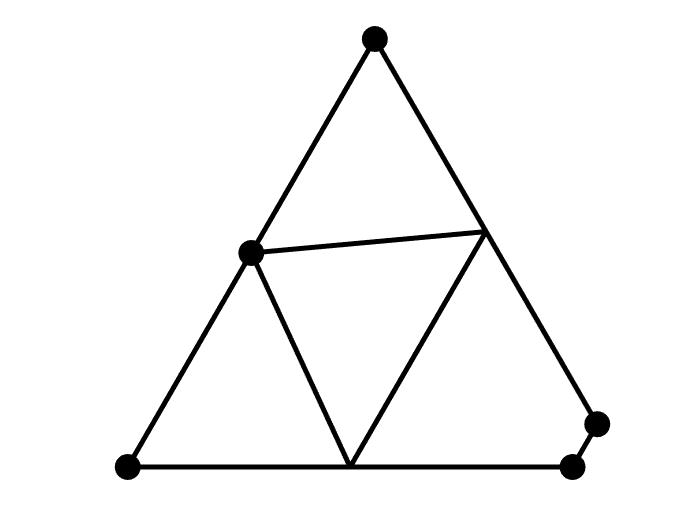}
    \includegraphics[width=0.3\linewidth]{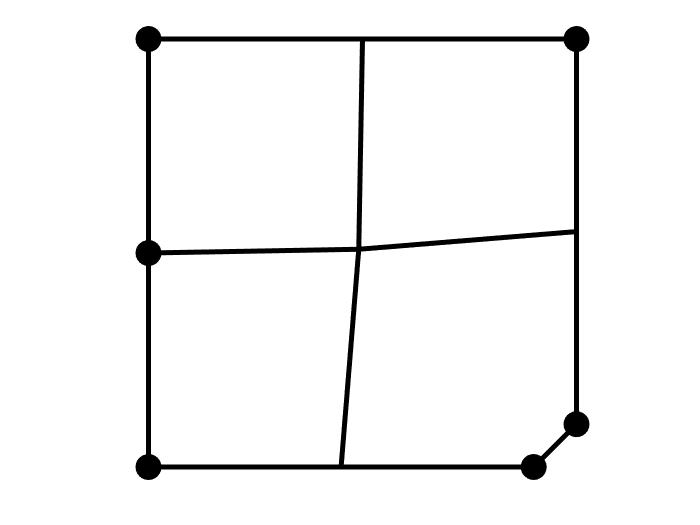}
    \caption{Top: the two polygons have been refined based on employing the ``plain" MP rules. Bottom: the two polygons have been first classified to belong to the class of ``triangular" and ``quadrilateral" element, respectively, and then refined accordingly. The vertices of the original polygons are marked with black dots.}
    \label{fig:midpoint-ANN}
\end{figure}

\begin{figure}%[!ht]
    \centering
    \includegraphics[width=\linewidth]{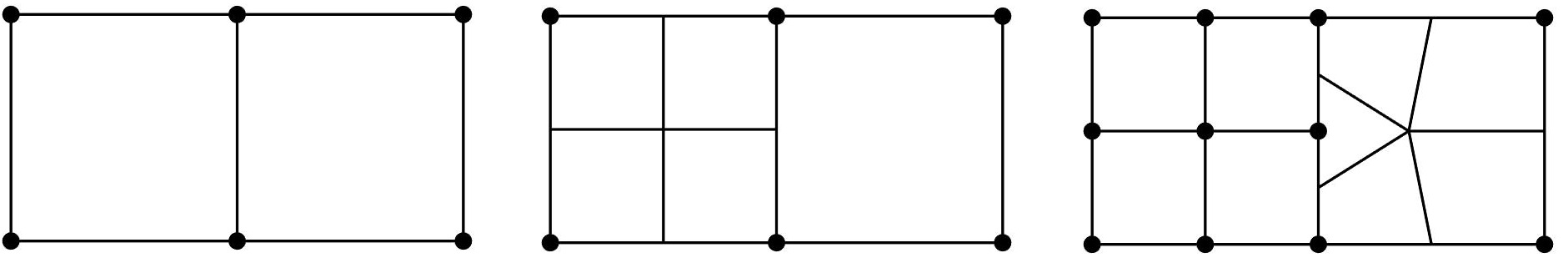}
    \caption[Midpoint strategy on adjacent elements.]{Left: the initial grid. Middle: the MP refinement strategy has been applied to the square on the left, therefore adding one node to the adjacent square on the right. Right: the MP refinement strategy has been applied also to the square on the right, dividing it into five sub-elements. The vertices of the grids that are employed to apply the MP refinement strategy of each stage are marked with black dots.}
    \label{fig:midpoint hanging nodes effect}
\end{figure}

In this section we discuss the problem of correctly identify the ``shape" of a general polygon, in order to later apply a suitable refinement strategy according to the chosen label of the classification. We start by observing that for polygons with ``regular" shapes, e.g. triangles, squares, pentagons, hexagons and so on, we can define ad-hoc refinement strategies. For example, satisfactory refinement strategies for triangular and quadrilateral elements can be designed in two dimensions, as shown in Figure~\ref{fig:regualar} (left). If the element $K$ is a triangle, the midpoint of each edge is connected to form four triangles; if $K$ is a square, the midpoint of each edge is connected with the centroid of the vertices of the polygon (to which we will refer as ``centroid", for short), i.e. the arithmetical average of the vertices coordinates, to form four squares.
For a regular polygon $K$ with more than fours edges, suitable refinement strategies can also be devised, see e.g. \cite{lai2016recursive} and Figure~\ref{fig:regualar} (right). The idea is to:
\begin{enumerate}
    \item Construct a suitably scaled and rotated polygon $\hat{K}$ with centroid that coincide with the centroid of the initial polygon $K$.
    \item Connect the vertices of $\hat{K}$ with the midpoints of the edges $K$, forming a pentagon for each vertex of the original polygon $K$.
\end{enumerate}
This strategy induces a partition with as many elements as the number of vertices of the original polygon $K$ plus one. The above refinement strategies for regular polygons have the following advantages:
\begin{itemize}
    \item They produce regular structures, thus preserving mesh quality.
    \item They enforce a modular structure, as the new elements have the same structure of the original one, easing future refinements.
    \item They keep mesh complexity low, as they add few vertices and edges.
    \item They are simple to be applied and have a low computational cost.
\end{itemize}
The problem of refining a general polygonal element is still subject to ongoing research. A possible strategy consist of dividing the polygon along a chosen direction into two sub-elements \cite{berrone2021refinement}; this strategy is very simple and has an affordable computational cost, and it also seems to be robust and to preserve elements' regularity. Because of its simplicity, however, this strategy can hardly exploit particular structures of the initial polygon.\\
\newline
Another possibility is to use a Voronoi tessellation \cite{hoshina2018simple}, where some points, called seeds, are chosen inside the polygon $K$ and each element of the new partition is the set of points which are closer to a specific seed, as shown in Figure~\ref{fig:voronoi ref}.
It is not obvious how many seeds to use and where to place them, but the resulting mesh elements are fairly rounded. The overall algorithm has a consistent but reasonable computational cost.\\
\newline
Another choice is to use the Mid-Point (MP) strategy, which consist in connecting the midpoint of each edge of $K$ with the centroid of $K$, as shown in Figure~\ref{fig:midpoint-ANN} (top).
This strategy is very simple and has a low computational cost. Modularity is enforced, in the sense that the resulting elements of the mesh are all quadrilaterals. %However, nodes are added to adjacent elements, as shown in Figure~\ref{fig:midpoint hanging nodes effect}, unless suitable (geometric) checks are included.
Notice that nodes are added to adjacent elements, therefore chaining the way those elements are refined unless suitable (geometric) checks are included, as shown in Figure~\ref{fig:midpoint hanging nodes effect}.
The main drawback of this strategy is that it potentially disrupts mesh regularity and the number of mesh elements increases very rapidly. Therefore, which refinement strategy is the most effective depends on the problem at hand and the stability properties of the numerical scheme used for its approximation.\\

\begin{figure}%[!ht]
    \centering
    \includegraphics[width=0.22\linewidth]{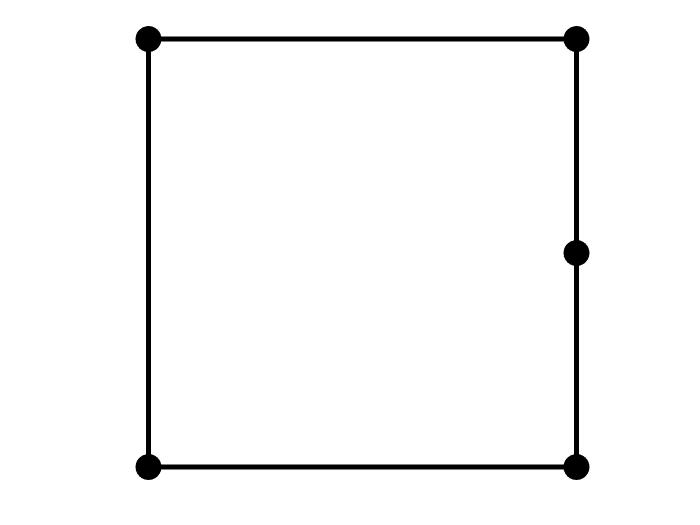}
    \includegraphics[width=0.22\linewidth]{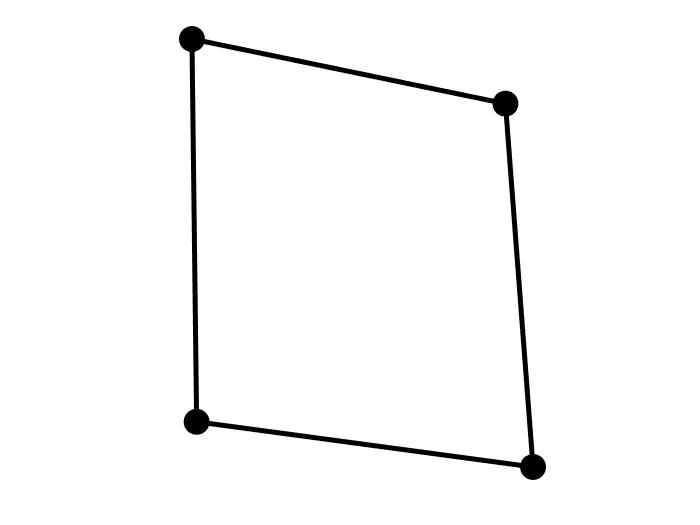}
    \includegraphics[width=0.22\linewidth]{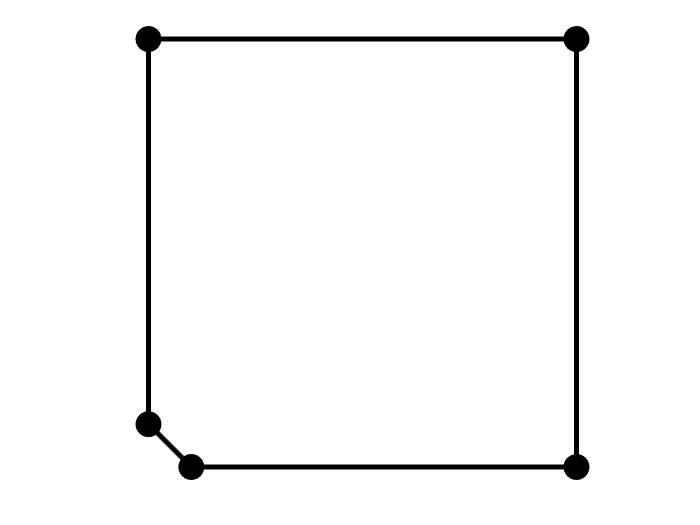}
    \includegraphics[width=0.22\linewidth]{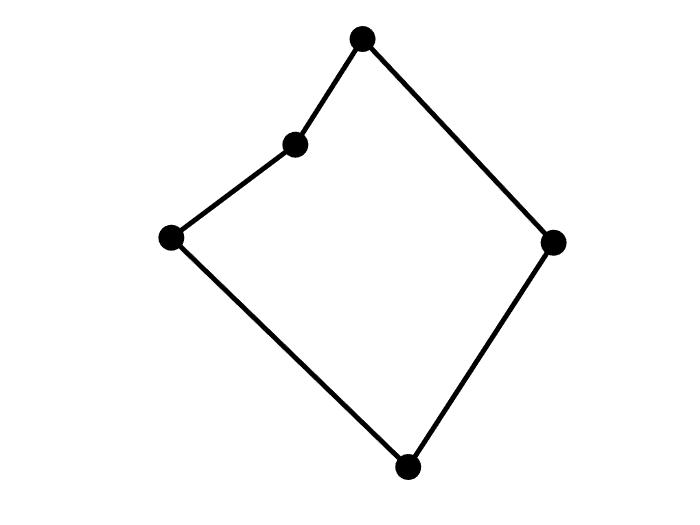}
    \caption[Polygons belonging to the class of squares.]{Illustrative examples of polygons belonging to the class of ``squares". They have been obtained by adding small distortions or extra aligned vertices to the reference square, plus rotations and scalings in some cases. The vertices are marked with black dots.}
    \label{fig:square class}
\end{figure}

\begin{figure}
\centering
\includegraphics[width=0.32\linewidth]{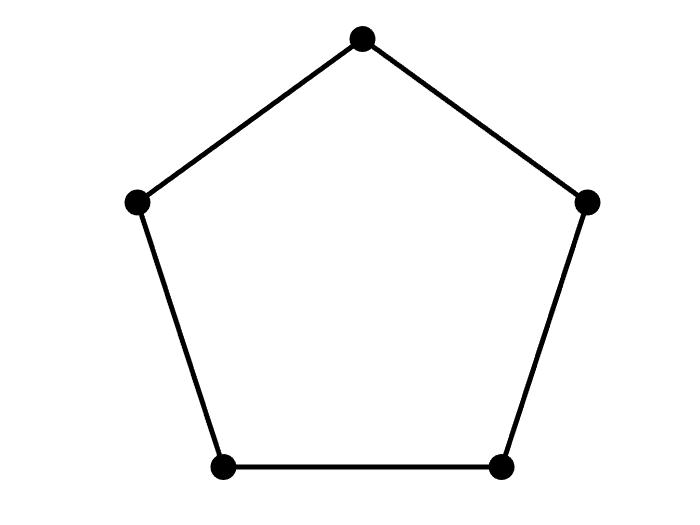}
    \includegraphics[width=0.25\linewidth]{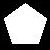}
    \caption[Binary image of a pentagon.]{Binary image of a pentagon of size {\color{black}$64 \times 64$ pixels}. Each pixel has value 1 (white pixel) if it is inside the polygon and 0 otherwise (black pixel). The binary images of each mesh element are then employed to classify the shape of the element, avoiding the automatic and exclusive use of geometric information.}
    \label{fig:binary image}
\end{figure}

\newline
In order to suitably drive these refinement criteria, we propose to use CNNs to predict to which ``class of equivalence" a given polygon belongs to. For example in Figure~\ref{fig:midpoint-ANN} (top) we show two elements refined using the MP strategy. They are clearly a quadrilateral and a pentagon, respectively, but their shapes are very similar to a triangle and a square respectively, and hence they should be refined as in Figure~\ref{fig:midpoint-ANN} (bottom). Loosely speaking, we are trying to access whether the shape of the given polygon is ``more similar" to a triangle, or a square, or a pentagon, and so on. The general algorithm is the following:
\begin{enumerate}
    \item Let $\mathcal{P} = \{P_1,\ P_2,\ ...\ P_m\}$ be a set of possible polygons to be refined and let $\mathcal{R} = \{R_1,\ R_2,\ ...\ R_n\}$ be a set of possible refinement strategies.
    \item We build a classifier $F: \mathcal{P} \rightarrow \mathcal{R}$ in such a way that suitable mesh quality metrics are preserved \cite{attene2019benchmark,VEMquality2021} (quality metrics will be described later in Section~3.3). The set of all polygons mapped into the same refinement strategy is a ``class of equivalence".
\end{enumerate}
In principle, any classifier $F$ may be used. However, understating the specific relevance of different geometric features of a general polygon (e.g. number of edges, area, etc...) might not be enough to operate a suitable classification. Instead, we can construct a database of polygons that can be used to ``train" our classifier $F$. In order to construct such database we proceed as follows. Starting from the ``reference" polygon in a class (e.g. the reference triangle for the class ``triangles", the reference square for the class ``squares" and so on) we generate a set of perturbed elements that still belong to the same class and are obtained by adding new vertices and/or adding noise to them, introducing rotations and so on. An illustrative example in the case of the class ``squares" is show in Figure~\ref{fig:square class}.
Training a function with labelled data is known as ``supervised learning", and CNNs are particularly well suited to deal with image classification problems in this context. Indeed, because of the inherently ``image classification" nature of this problem, polygons can be easily converted into binary images: each pixel assumes value 1 if it is inside the polygon and 0 otherwise, as shown in Figure~\ref{fig:binary image}.

{\color{black}A binary image of a polygon can be efficiently generated as follows:
		\begin{enumerate}
			\item Scale and translate the original polygon into the reference box $(0,1)^2$.
			\item Construct a grid with the desired number of pixels.
			\item From each edge of the polygon sample properly spaced points and assign value 1 to pixel containing those points. This will define an inner and an outer region of pixels.
			\item Assign value 1 to each pixels in the inner region: starting from an inner pixel, recursively assign value 1 to neighbouring pixels until other pixels with value 1 are met.
		\end{enumerate}
	For the sake of computational efficiency, this process may also be performed approximately to a certain extent. Indeed, the neural network directly learns from the image representation of the polygon and, if trained properly, it should be robust with respect to small distortions.

% For the sake of computational complexity, this operation may also be performed approximately, since the CNN needs to be robust with respect to image noise.
Identifying the ``shape" of a polygonal element based on its geometrical properties is also possible using the so-called Shape Analysis (see e.g. \cite{dryden2016statistical}). This approach relies on rotating and translating the original element so as to minimize a suitable distance function, called ``Procrustes distance", from a reference shape. The reference shape with the smallest Procrustes distance is then selected. This is a viable option in two dimensions, because the vertices of a polygon can be naturally ordered clock-wise, whereas in three dimensions it has much higher computational cost. Despite the fact that different classification algorithms may be more effective for different situations, depending on the requirements on accuracy, generalization and computational costs, we decided to use a machine learning based approach because:
		\begin{itemize}
			\item The behaviour of the classifier may be tuned using examples. Therefore, an explicit description of the reference shape is not required. This allows to define very general "classes of equivalence". For example, an alternative database may be generated as follows:
			\begin{itemize}
				\item pick a polygon and apply different refinement strategies
				\item rank the outcomes according to some criterion;
				\item assign the label corresponding to the refinement strategy which attained the highest score.
			\end{itemize}		
			\item If a new set of elements is not classified as desired, a simple solution is to label thoses samples and include them in the training process of the CNN. Solving this problem using Shape Analysis is more complex, because the set of possible reference shapes need to be re-design.
			\item Labeling a new sample using a neural network is computationally efficient, especially when dealing with multiple classes: the CNN extracts key features from the image, which are then used to assign the label. Instead, using a Shape Analysis approach the new sample is compared with every reference shape separately."
		
		\end{itemize}
}

\subsection{Supervised learning for image classification}
Consider a two dimensional gray-scale image, represented by a tensor $B \in \mathbb{R}^{m \times n},\ m,n\geq 1,$ and the corresponding label vector $y = ([y]_j)_{j =1,..,\ell} \in  [0,1]^\ell,$
where $\ell \geq 2$ is the total number of classes, and $[y]_j$ is the probability of $B$ to belong to the class $j$ for $j = 1:\ell$. For the case of polygons classification $B \in \{0,1\}^{m \times n}$, i.e images are binary. Moreover, in our case the classes are given by the label ``triangle", ``square", ``pentagon" and so on.\\
In a supervised learning framework, we are given a dataset of desired input-output couples $\{(B_i,y_i)\}_{i=1}^N$, where $N$ is number of labelled data.
We consider then an image classifier represented by a function of the from $F: \mathbb{R}^{m \times n} \rightarrow~(0,1)^\ell,$ in our case a CNN, parameterized by $w \in \mathbb{R}^M$ where $M\geq 1$ is the number of parameters. Our goal is to tune $w$ so that $F$ minimizes the data misfit, i.e.
$$\min_{w \in \mathbb{R}^M} \sum_{i \in I} l(F(B_i),y_i),$$   
where $I$ is a subset of $\{1,2,...,N\}$ and $l$ is the cross-entropy loss function defined as
$$ l(F(B),y) = \sum_{j = 1}^\ell -[y]_j\log[F(B)]_j.$$
This optimization phase is also called ``learning" or ``training phase". During this phase, a known shortcoming is ``overfitting": the model fits very well the data used in the training phase, but performs poorly on new data. For this reason, the data set is usually splitted into: i) training set: used to tune the parameters during the training phase; % In our case, we used 60\% of the whole data set.
ii) validation set: used to monitor the model capabilities on different data during the training phase. The training is halted if the error on the validation set starts to increase; % In our case, we used 20\% of the whole data set.
iii) test set: used to access the actual model performance on new data after the training. % In our case, we used the remaining 20\% of the data.
\\
While the training phase can be computationally demanding, because of the large amount of data and parameters to tune, it needs to be performed off line once and for all. Instead, classifying a new image using a pre-trained model is computationally fast: it requires only to evaluate $F$ on a new input. The predicted label is the one with the highest estimated probability.

\begin{figure}
    \hspace{-2.2cm}
    \includegraphics[width = 1.4\linewidth]{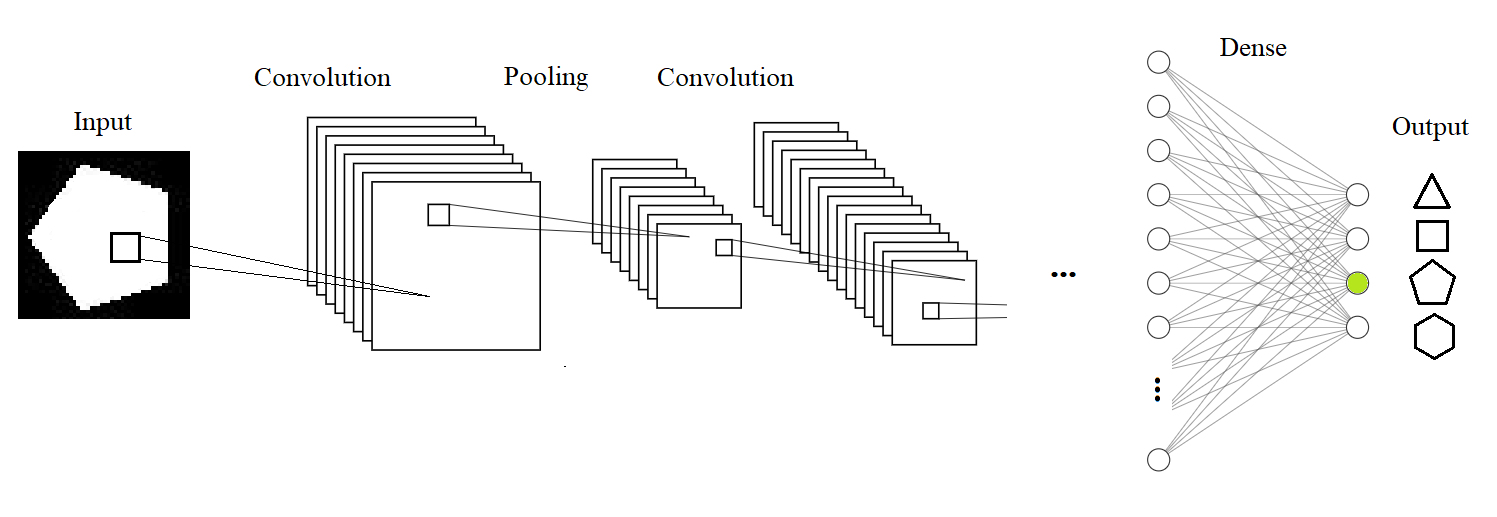}
    \caption{Simplified scheme of a CNN architecture employed for classification of the shape of polygons. The large squares represent the image channels after applying a convolution layer or a pooling layer, while the small squares represent layer filters scanning through every channel. After applying a convolution layer the number of channels is multiplied by the number of  feature maps, while after applying a pooling layer the size of each channel diminishes accordingly. Circles represent neurons, one for each input and output, and connections represent linear dependencies. Finally, labels are represented using geometrical shapes.}
    \label{fig:CNN architecture}
\end{figure}

{\color{black}
\subsection{Convolutional Neural Networks}
CNNs are parameterized functions, in our case of the form $$\text{CNN}:~\mathbb{R}^{m \times n}~\rightarrow (0,1)^\ell,\hspace{0.5cm} m,n,\ell \geq 1,$$ constructed by composition of simpler functions called ``layers of neurons". We are now going to introduce different types of layers:
\begin{itemize}

    \item Convolutional layers: linear mappings of the form $\textsc{Conv}: \mathbb{R}^{m \times n\times c} \rightarrow \mathbb{R}^{m \times n\times \Bar{h}}$ with $m,n,c,\Bar{h}\geq 1$, where $m$ and $n$ are the size of the input image, $c$ is the number of channels, e.g. $c = 3$ for a colored image, and $\Bar{h}$ is the number of features maps. For an image $B \in \mathbb{R}^{m \times n}$ and a kernel $K \in \mathbb{R}^{(2k+1) \times (2k+1)}$ the convolution operator $*$ is defined as
    $$ [K * B]_{i, j}=\sum_{p, q=-k}^{k} [K]_{k+1+p, k+1+q} [B]_{i+p, j+q}, \quad i=1:m,\ j=1:n,$$
    with zero padding, i.e. $B_{i+p, j+q} = 0$ when indexes are out of range. This operation can be viewed as a filter scanning through image $B$, extracting local features that depend only on small subregions of the image. This is effective because a key property of images is that close pixels are more strongly correlated than distant ones. The scanning filter mechanism provides the basis for the invariance of the output to translations and distortions of the input image \cite{bishop2006pattern}. The convolutional layer is defined as
    $$[\textsc{Conv}(B)]_{i}=\sum_{j=1}^{c} [K]_{:,:,i,j} *[B]_{:,:,j} + [b]_{i} \mathbf{1}, \quad i=1: \Bar{h},$$
    where the colon index denotes that all the indexes along that dimension are considered, $\mathbf{1} \in \mathbb{R}^{m \times n}$ is the $m \times n$ matrix with all entries equal to 1, $K \in \mathbb{R}^{(2k+1) \times (2k+1) \times \Bar{h} \times c}$ is a kernel matrix and $b \in \mathbb{R}^{\Bar{h}}$ is a bias vector of coefficients to be tuned.

    \item Batch normalization ($\textsc{Norm}$): linear mappings used to speed up training and reduce the sensitivity to network initialization \cite{ioffe2015batch}.

    \item Pooling layers: mappings used to perform down-sampling, such as $$\textsc{Pool}: \mathbb{R}^{m \times n \times c} \rightarrow \mathbb{R}^{\lceil\frac{m}{s}\rceil \times \lceil\frac{n}{s}\rceil \times c} \hspace{0.5cm}m,n,c \geq 1,$$ $$[\textsc{Pool}(B)]_{i,j,t}=\max_{p, q=1:k} B_{s(i-1)+p, s(j-1)+q,t}$$ with $s \geq 1$ and zero padding. They improve the invariance of the output with respect to translations of the input.
    
    \item Activation functions: mappings used to introduce non-linearity, such as the rectified linear unit $[\textsc{RELU}(x)]_i = \max(0,x_i)$.

    \item Dense layers: generic linear mappings of the form $\textsc{Linear}:  \mathbb{R}^{m \times n\times c} \rightarrow \mathbb{R}^\ell,$ $n,m,c,\ell \geq 1$ defined by parameters to be tuned. They are used to separate image features extracted in the previous layers.
    
    \item $\textsc{Softmax}: \mathbb{R}^\ell \rightarrow (0,1)^\ell,$ where $\ell\geq 1$ is the number output classes, $[\textsc{Softmax}(x)]_i = \frac{e^{x_i}}{\sum_{j=1}^\ell e^{x_j}}$. They are used to assign a probability to each class.
   
\end{itemize}
In practise, subsequent application of convolutional, activation and pooling layers may be used to obtain a larger degree of invariance to input transformations such as rotations, distortions, etc. A visual representation of a CNN is shown in Figure~\ref{fig:CNN architecture} for the case of regular polygons classification.}

\section{CNN-enhanced refinement strategies}
In this section we present two strategies to refine a general polygon that exploit a pre-classification step of the polygon shape. More specifically, we assume that a CNN for automatic classification of the polygon label is given. The first strategy consists in enhancing the classical MP algorithm, whereas the second strategy exploits the refinement criteria for regular polygons illustrated in Section 2.
%%%%%%%%%%%%%%%%%%%%%%%%%%%%%%%%%%%%%%%%%%%%%%%%%%%%%%%%%%%%%%
\subsection{A CNN-enhanced MP strategy}
Assume we are given a general polygon $P$ to be refined  and its label $\mathcal{L}$, obtained using a CNN for classification of polygon shapes. Here $\mathcal{L} \geq 3$ is an integer, where $\mathcal{L}=3$ corresponds to the label ``triangle", $\mathcal{L}=4$ corresponds to the label ``square", and so on.
\begin{table}
        % \vspace{-3.2cm}
        % \hspace{-3cm}
        \begin{tabular}{M{27mm}M{27mm}M{27mm}M{27mm}}
        $\mathcal{L} = 3$ (triangle) & $\mathcal{L} = 4$ (quad) & $\mathcal{L} = 5$ (pentagon) & $\mathcal{L} = 6$ (hexagon) \\

    \includegraphics[width=\linewidth]{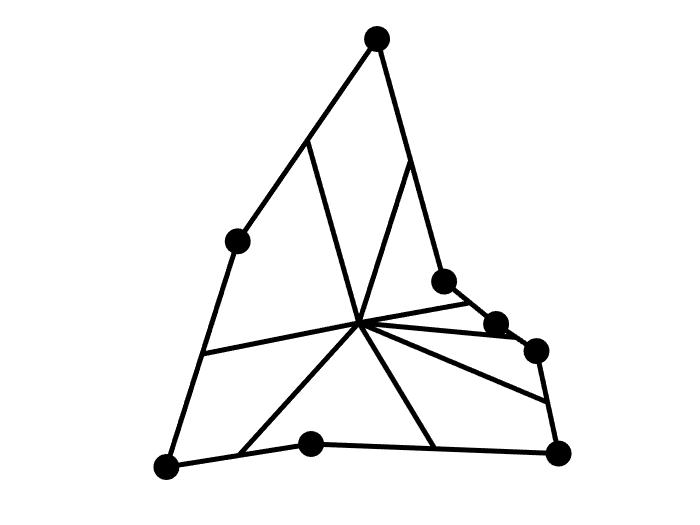}&
    \includegraphics[width=\linewidth]{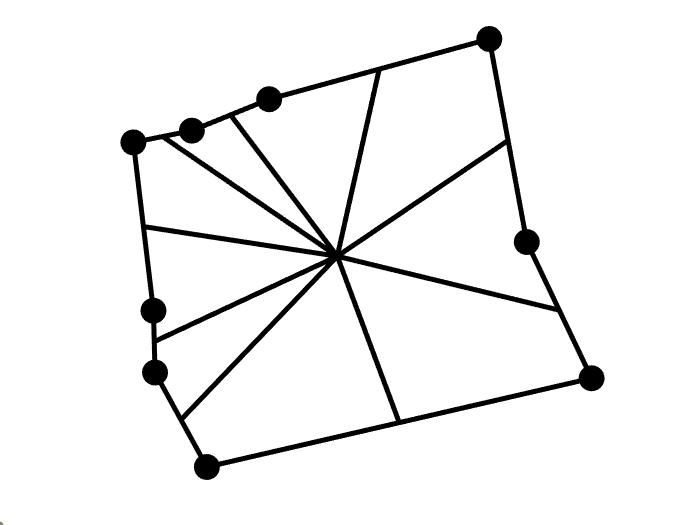}&
    \includegraphics[width=\linewidth]{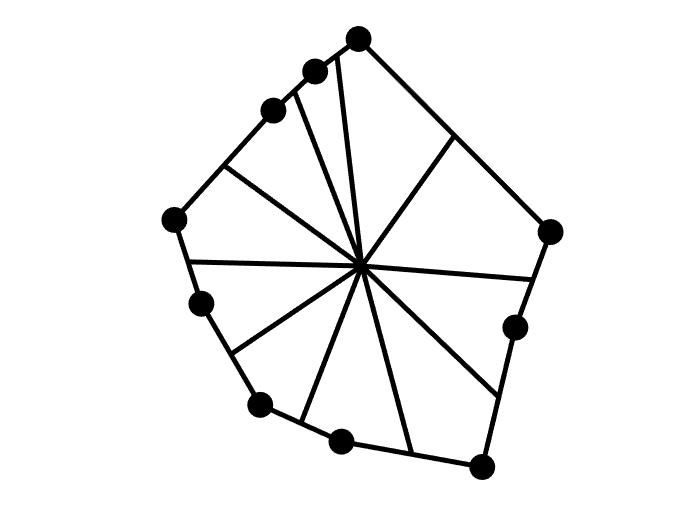}&
    \includegraphics[width=\linewidth]{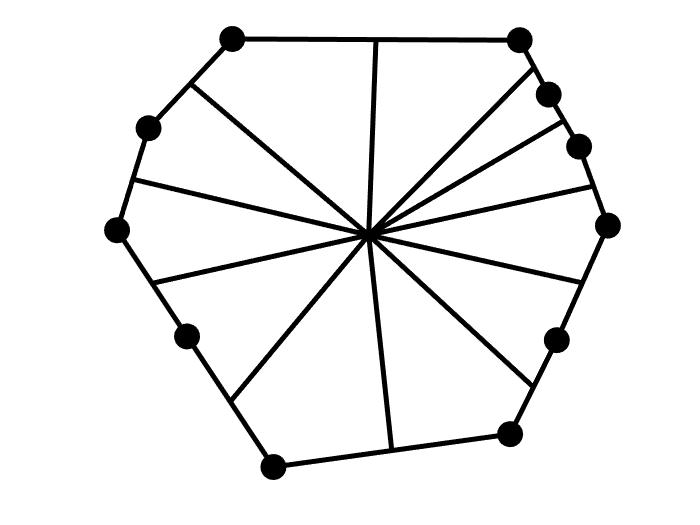}\\
        
         \includegraphics[width=\linewidth]{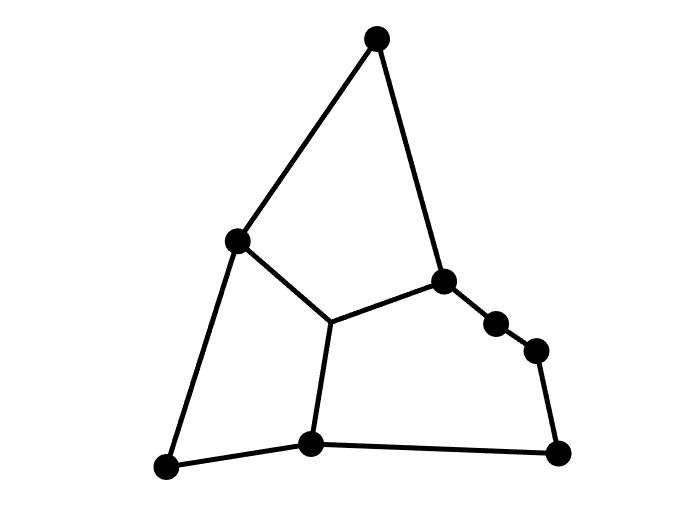}&
    \includegraphics[width=\linewidth]{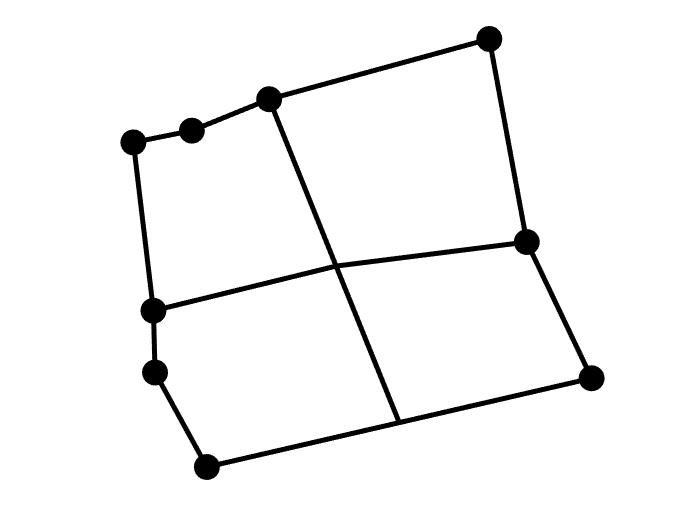}&
    \includegraphics[width=\linewidth]{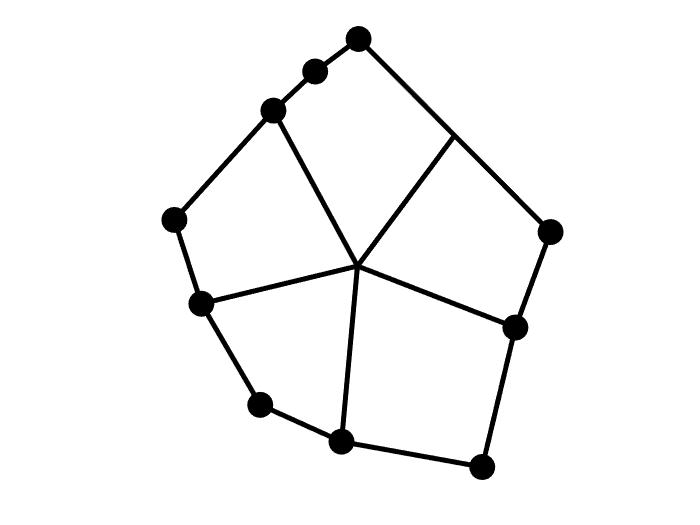}&
    \includegraphics[width=\linewidth]{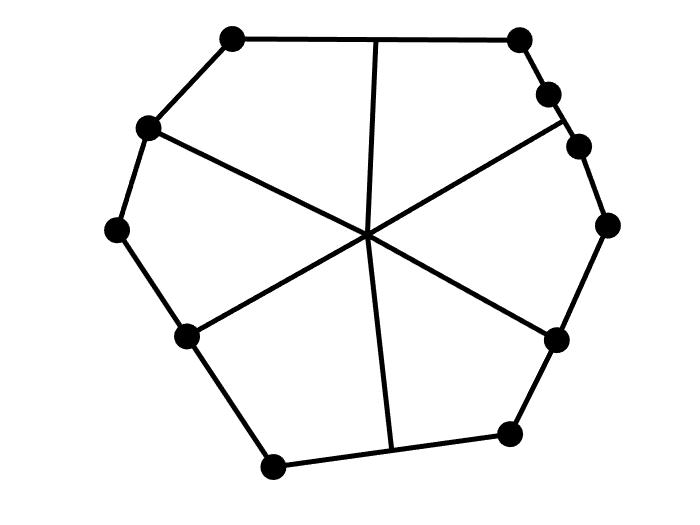}\\
    
    \includegraphics[width=\linewidth]{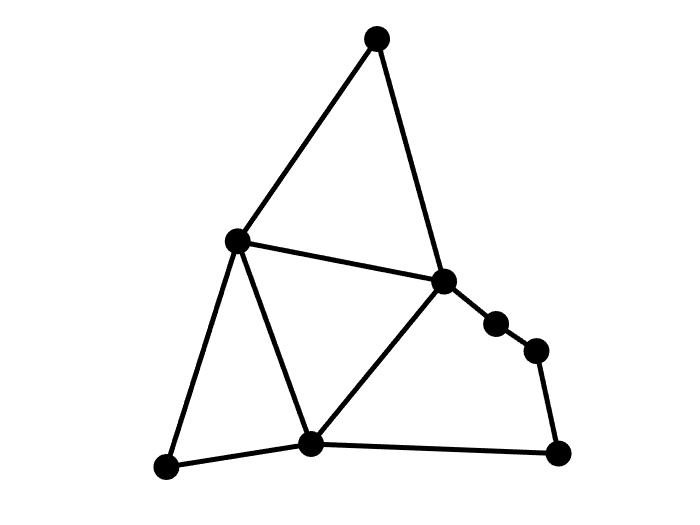}&
    \includegraphics[width=\linewidth]{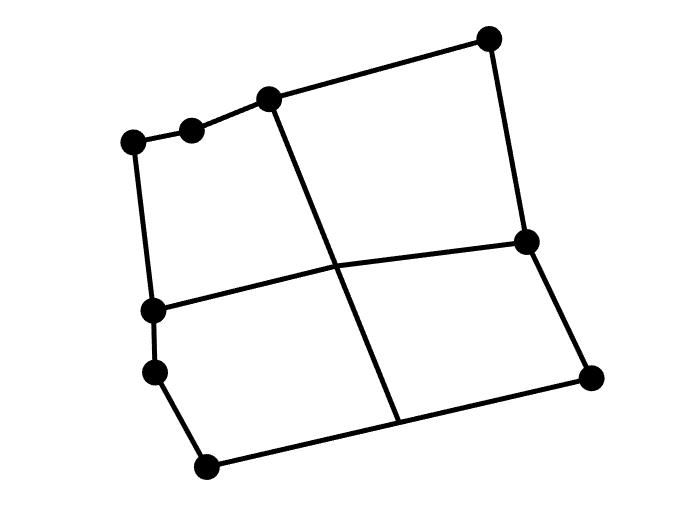}&
    \includegraphics[width=\linewidth]{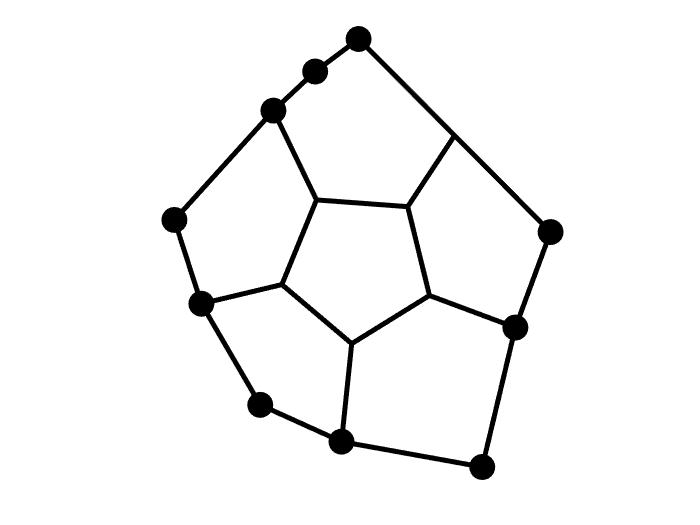}&
    \includegraphics[width=\linewidth]{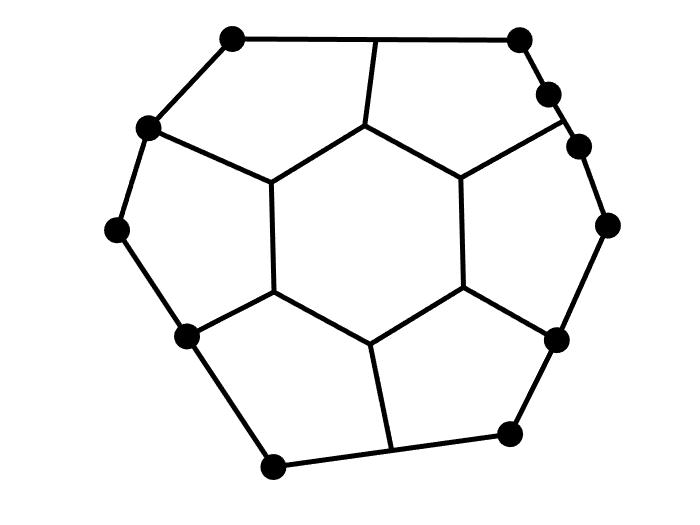}\\

        \end{tabular}
        
                %\captionsetup{width=1.25\linewidth}
        \captionof{figure}{Samples of polygons refined using the MP (top), CNN-MP (middle) and CNN-RP (bottom) refinement strategies. The elements have been classified using the CNN algorithm with labels $\mathcal{L} = 3,4,5,6$, respectively. 
    The vertices of the original polygons are marked with black dots.}
       \label{fig:strategies}
\end{table}
If the polygon $P$ has a large number of (possibly aligned) vertices, applying the MP strategy may lead to a rapid deterioration of the shape regularity of the refined elements. In order to reduce the number of elements produced via refinement and to improve their quality, a possible strategy is to enhance via CNNs the MP refinement strategy, and apply the MP refinement strategy not to the original polygon $P$ but to a suitable approximate polygon $\hat{P}$ with a number of vertices $\mathcal{L}$ identified by the CNN classification algorithm, as described in Algorithms 1 and 2.
\begin{algorithm}[t]
\SetAlgoLined
\LinesNumbered
\KwIn{polygon $P$.}
\KwOut{partition of $P$ into polygonal sub-elements.}
\DontPrintSemicolon
\BlankLine
\BlankLine
Convert $P$, after a proper scaling, to a binary image as shown in Figure~\ref{fig:binary image}. \;
Apply a CNN for classification of the polygon shape and obtain its label $\mathcal{L} \geq 3$. \;
Based on $\mathcal{L}$, identify the refinement points on the boundary of $P$, as described in Algorithm 2. \;
Connect the refinement points of $P$ to its centroid $c_P$. \;
\caption{CNN-enhanced Mid-Point (CNN-MP) refinement strategy}
\end{algorithm}
\begin{algorithm}[t]
\SetAlgoLined
\LinesNumbered
\KwIn{polygon $P$, label $\mathcal{L} \geq 3$.}
\KwOut{refinement points on the boundary of $P$.}
\DontPrintSemicolon
\BlankLine
\BlankLine
Build a polygon $\hat{P}$ suitably approximating $P$: select $\mathcal{L}$ vertices $\hat{v}_1,\hat{v}_2, ... \hat{v}_\mathcal{L}$, among the vertices of $P$, that maximize $\sum_{i,j = 1}^\mathcal{L} ||\hat{v}_i-\hat{v}_j||$.\;
Compute the centroid $c_P$ of $P$, the edge midpoints of $P$ and the edge midpoints $\{\hat{m}_i\}_{i = 1}^\mathcal{L}$ of $\hat{P}$.\;
For every edge midpoint $\hat{m}_i$ of $\hat{P}$ : find the closest point to $\hat{m}_i$ and $c_P$, among the vertices and the edges midpoints of $P$.\;
\caption{Identification of the refinement points}
\end{algorithm}
Examples of refined polygons using the MP strategy are shown in Figure~\ref{fig:strategies} (top) whereas the analogous ones obtained employing the CNN-enhanced Mid-Point (CNN-MP) refinement strategy are shown in Figure~\ref{fig:strategies} (middle). We point out that the computational cost of the CNN-MP strategy is very low. Moreover, parallelism is enforced in a stronger sense: the CNN does not distinguish between one edge of a polygon and two aligned edges, solving the problem of refining adjacent elements pointed out in Figure~\ref{fig:midpoint hanging nodes effect}.

%%%%%%%%%%%%%%%%%%%%%%%%%%%%%%%%%%%%%%%%%%%%%%%%%%%%%%%%%%%%%%

\subsection{A new ``reference polygon'' based refinement strategy}
Assume, as before, that we are given a general polygon $P$ to be refined together with its label $\mathcal{L}\geq 3$, that can be obtained employing a CNN classification algorithm.
If the given polygon is a reference polygon we could refine it based on employing the refinement strategies described in Section 2 and illustrated in Figure~\ref{fig:regualar}, where the cases $\mathcal{L} = 3,4,5,6$ are reported. Our goal is to extend these strategies so that they can be applied to general polygons. In order to do that, the idea of the algorithm is to first compute the refinement points of $P$, as before, and then connect them using the refinement strategy for the class $\mathcal{L}$. More precisely, our new CNN-enhanced Reference Polygon (CNN-RP) refinement strategy is described in Algorithms 2 and 3 and illustrated in Figure~\ref{fig:strategies} (bottom). Notice that lines 1-2-3 in Algorithm 3 are the same as in Algorithm 1.
%%%%%%
\begin{algorithm}[th!]
\SetAlgoLined
\LinesNumbered
\KwIn{polygon $P$.}
\KwOut{partition of $P$ into polygonal sub-elements.}
\DontPrintSemicolon
\BlankLine
\BlankLine
Convert $P$, after a proper scaling, to a binary image as shown in Figure~\ref{fig:binary image}. \;
Apply a CNN for classification of the polygon shape and obtain its label $\mathcal{L} \geq 3$. \;
Based on $\mathcal{L}$, identify the refinement points on the boundary of $P$, as described in Algorithm 2. \;
\If{$\mathcal{L} = 3$}{
 Connect the  refinement points of $P$ so as to form triangular sub-elements.\;
}
\BlankLine
\If{$\mathcal{L} = 4$}{
Connect the refinement points of $P$ to its centroid, so as to form quadrilateral sub-elements. \;
}
\BlankLine
\If{$\mathcal{L} \geq 5$}{
Construct inside of $P$ a suitably scaled and rotated regular polygon with $\mathcal{L}$ vertices and with the same centroid of $P$.\;
Connect the vertices of the inner regular polygon with the refinement points of the outer polygon $P$, so as to form sub-elements as shown in Figure~\ref{fig:regualar}.\;
}
\caption{CNN-enhanced Reference Polygon (CNN-RP) refinement strategy}
\end{algorithm}
This strategy can be applied off line with a low computational cost and has the advantage to enforce parallelism as each mesh element can be refined independently.\\

%%%%%%%%%%%%%%%%%%%%%%%%%%%%%%%%%%%%%%%%%%%%%%%%%%%%%%%%%%%%%

Notice that for a non-convex polygon, the CNN-RP and the CNN-MP strategies do not guarantee in general to generate a valid refined element, because the centroid could lie outside the polygon. In practise, they work well even if the polygon is ``slightly" non-convex. However, in case of a non-valid refinement, it is always possible to first subdivide the polygon into two elements, possibly of comparable size, by connecting two of its vertices. In general, non-valid refinements may be detected because they are not valid partitions, i.e. when the new elements are overlapping or when they do not cover the original element. In practise, one may simply check whether or not some points lie inside the original polygon and whether or not edges are intersecting each other.

\subsection{Quality metrics}
\begin{figure}
    \centering
    \includegraphics[width = 0.35\linewidth]{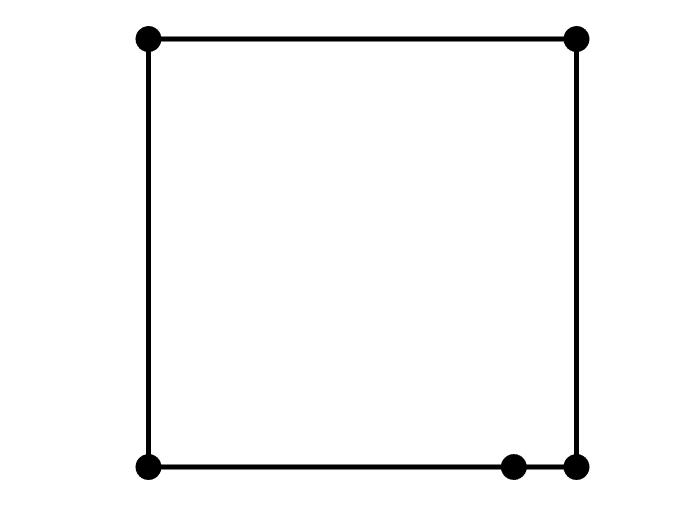}
    \caption[A regular polygon with small Edge Ratio and Normalized Point Distance.]{A polygon with a small edge. Although its shape is regular, ER and NPD metrics assume small values, depending on the size of the smallest edge.}
    \label{small ER}
\end{figure}

In order to evaluate the quality of the proposed refinement strategies, we recall some of the quality metrics introduced in \cite{attene2019benchmark}. The diameter of a domain $\mathcal{D}$ is defined, as usual, as $\textrm{diam}(\mathcal{D}) := \sup \{|x-y|,\ x,y \in \mathcal{D}\}.$ Given a polygonal mesh, i.e. a set of non-overlapping polygonal regions $\{P_i\}_{i=1}^{N_P}$, $N_P \geq 1$, that covers a domain $\Omega$, we can define the mesh size $h = \max_{i=1:N_P} \textrm{diam}(P_i).$ For a mesh element $P_i$, the Uniformity Factor (UF) is defined as $\textrm{UF}_i = \frac{\textrm{diam}(P_i)}{h}$, $i = 1,...,N_P$. This metric takes values in $[0,1]$.\\
\newline
For a polygon $P$, we also introduce the following quality metrics, taken from \cite{attene2019benchmark}:
\begin{enumerate}
    %\item Uniformity Factor (UF): $\frac{\textrm{diam}(P)}{h}$.
    \item Circle Ratio (CR): ratio between the radius of the inscribed circle and the radius of the circumscribed circle of $P$:
    $$\frac{\max_{\{B(r) \subset P\}}  r }{\min_{\{P \subset B(r)\}} r },$$
	where $B(r)$ is a ball of radius $r$. For the practical purpose of measuring the roundness of an element the radius of the circumscribed circle has been approximated with $\textrm{diam}(P) / 2$.
    \item Area-Perimeter Ratio (APR):
    \begin{equation*}
        \frac{4 \pi\ \textrm{area}(P)}{\textrm{perimeter}(P)^2}.
    \end{equation*}
    
    \item Minmum Angle (MA): minimum inner angle of $P$, normalized by $180^\circ$.
    \item Edge Ratio (ER): ratio between the shortest and the longest edge of $P$. 
    \item Normalized Point Distance (NPD): minimum distance between any two vertices of $P$, divided by the diameter of the circumscribed circle of $P$.
\end{enumerate}
These metrics are scale-independent and take values in $[0,1]$. The more regular the polygons are, the larger CR, APR and MA are. Large values of ER and NPD also indicate that the polygon is well proportioned and not skewed. However, small values of ER and NPD do not necessarily mean that the element is not shaped-regular, as shown in Figure~\ref{small ER}.\\
Notice that some quality metrics are not as important as others. For example, for several polytopal methods it is known that a small Edge Ratio does not necessarily deteriorates the
accuracy of the method \cite{brenner2018virtual,beirao2017stability,droniou2021robust,mu2015shape}. 

\section{CNN training}
\begin{figure}
    %\centering
    \hspace{-2.4cm}
    \includegraphics[width = 0.7\linewidth]{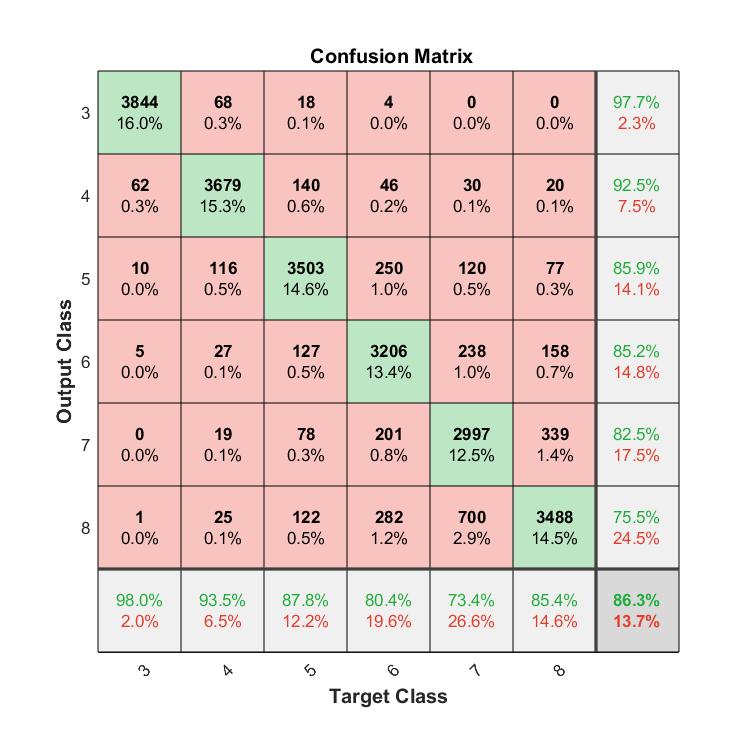}
    \includegraphics[width = 0.7\linewidth]{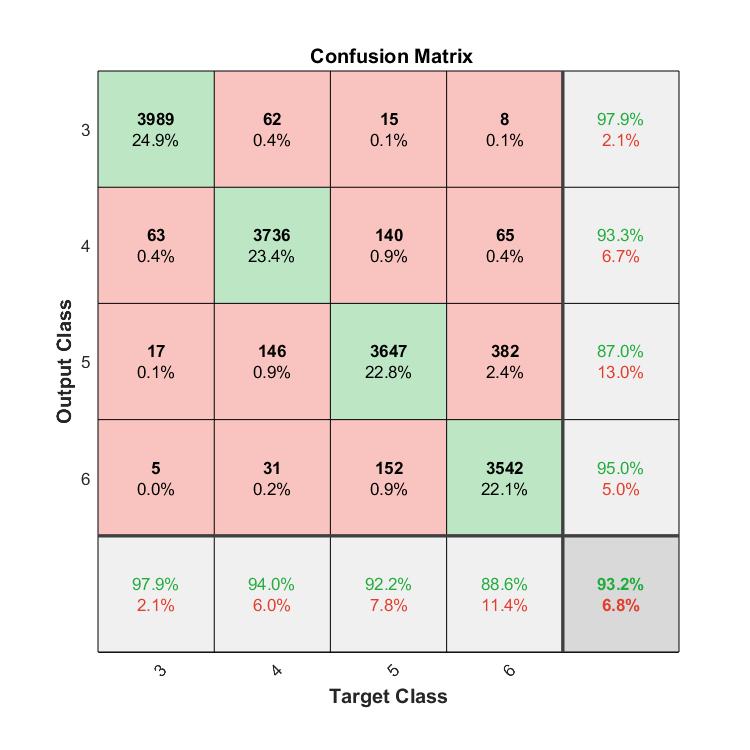}
    \caption{Confusion matrices for polygons classification.
    Left: the number of classes is $\ell = 6$  and the target classes vary from $\mathcal{L} = 3$ (triangles) to $\mathcal{L} = 8$ (octagons). Right: the number of classes is $\ell = 4$  and the target classes vary from $\mathcal{L} = 3$ (triangles) to $\mathcal{L} = 6$ (hexagons). The prediction accuracy for each target class decreases as more target classes are considered.}
    \label{fig:confusion}
    
\end{figure}

The CNN architecture we used for polygons classification is given by \\
{\color{black}
$ \text{CNN}: \{0,1\}^{64 \times 64} \rightarrow (0,1)^\ell,$
where
\begin{align*} 
\text{CNN} =\ &  \textsc{Conv}(k = 1, \Bar{h} =2)\  \rightarrow \textsc{Norm} \rightarrow \text{ReLU} \rightarrow  \textsc{Pool}(k = 2, s = 2)\rightarrow \\ 
& \textsc{Conv}(k = 1,\Bar{h} =4) \rightarrow \textsc{Norm} \rightarrow \textsc{ReLU}\rightarrow \textsc{Pool}(k = 2, s = 2) \rightarrow \\
& \textsc{Conv}(k = 1,\Bar{h} =8) \rightarrow \textsc{Norm} \rightarrow \textsc{ReLU}\rightarrow \textsc{Pool}(k = 2, s = 2) \rightarrow \\
& \textsc{Conv}(k = 1,\Bar{h} =16) \rightarrow \textsc{Norm} \rightarrow \textsc{ReLU}\rightarrow \textsc{Pool}(k = 2, s = 2) \rightarrow \\
& \textsc{Conv}(k = 1,\Bar{h} =32) \rightarrow \textsc{Norm} \rightarrow \textsc{ReLU}\rightarrow \textsc{Pool}(k = 2, s = 2) \rightarrow \\
& \textsc{Conv}(k = 1,\Bar{h} =64) \rightarrow \textsc{Norm} \rightarrow \textsc{ReLU}\rightarrow \textsc{Linear} \rightarrow \textsc{Softmax},
\end{align*}
}
where $k,\Bar{h},s$ are defined as in Section 2.2. \textcolor{black}{This choice of a deep architecture is motivated by the fact that many convolutional and pooling layers are needed to enforce the invariance to any possible rotation of the input polygon.} \textcolor{black}{The size of the images was selected to be $64\times 64$ pixels, i.e. a resolution large enough to apply five pooling layers, whose effect is to down-sample the image.} {\color{black}It was empirically observed that a larger resolution was not needed for the datasets we are going to consider. This is motivated by the fact that the approximation properties of neural networks generally depend on the dimension of the manifold where data are, and not on the dimension of the input space (see e.g. \cite{petersen2020neural}). A larger resolution maybe needed in order to be able to differentiate between polygons characterized by smaller angles, e.g. to correctly distinguish between a regular polygon with 9 sides an one with 10 sides.}\\
\textcolor{black}{For each class, we generated 20.000 images transforming regular polygons by adding edges and noise to the vertices. Since we will consider datasets with $\ell = 4$ and $\ell = 6$ classes, the total number of images will be respectively 80.000 and 120.000.} We set training, validation and test sets equal to 60\%-20\%-20\% of the whole dataset, respectively. \textcolor{black}{To train the neural networks we used the Adam (adaptive moment estimation) optimizer, see e.g. \cite{kingma2014adam}.}\\
Initially we selected the number of target classes to be equal to $\ell = 6$, i.e. polygons are sampled from triangles to octagons. We show the confusion matrix in Figure~\ref{fig:confusion} (left).
The same results obtained with $\ell = 4$, i.e. target classes varying from $\mathcal{L} = 3$ (triangles) to $\mathcal{L} = 6$ (hexagons), are shown in Figure~\ref{fig:confusion} (right). From these results it seems that the prediction accuracy is better in the case of a smaller set of target classes. This is expected, as for example a regular octagon is much more similar, in terms of angles amplitude and edges length, to a regular heptagon than to a regular triangle. Moreover, for polygons with many edges more pixels might be required in order  to appreciate the differences between them. In the following numerical experiments we have decided to choose $\ell = 4$, as this choice seems to balance the effectiveness of our classification algorithm with the computational cost. We also remark that for the following reasons:\\
\begin{itemize}
    \item Refining heptagons and octagons as if they were hexagons does not seem to affect dramatically the quality of the refinement.
    \item Ad-hoc refinement strategies for polygons with many edges seem to be less effective because more sub-elements are produced.
    \item A considerable additional computational effort might be required to include more classes.
    \item The more classes we use, the easier the possibility of a misclassification error is and hence to end up with a less robust refinement procedure.
\end{itemize}
\textcolor{black}{Considering polygon classes ranging from triangles to hexagons yields a satisfactory accuracy of 93.2\% as shown in Figure~\ref{fig:confusion} (right). Consider also that an accuracy close to 100\% is not realistic because when distorting for example a pentagon by adding noise to its vertices there is a considerable probability to turn it into something indistinguishable from a slightly distorted “square", which therefore we would like to classify as “square" even though the original label is “pentagon".\\
“Hexagons" and “pentagons" are the most misclassified elements, when classified as “pentagons" and “squares" by the CNN, respectively. In both cases, applying the wrong refinement strategy should not yield dramatic differences in terms of elements quality, as for example it would instead misclassifying an “hexagon" as a “triangle". Moreover, an effective refinement algorithm should be robust to some extent with respect to misclassification errors. These remarks will be confirmed by numerical experiments, shown in Section~6.}\\
Thanks to the limited number of dataset samples and network parameters, the whole algorithm (dataset generation, CNN training and testing) took approximately {\color{black} six minutes} using MATLAB2019b on a Windows OS 10 Pro 64-bit, Intel(R) Core(TM) i7-8750H CPU (2.20GHz / 2.21GHz) and 16 GB RAM memory.
\\ {\color{black}Improving the accuracy of the CNN using more data and/or a network with more layers is possible. However, the dataset will probably not be fully representative of mesh elements on which the CNN is going to be applied, because it is artificially generated. Therefore, attempting to achieve a very high accuracy may cause the CNN to overfit such dataset and hence to generalize poorly on real mesh elements. A different approach would be using directly real mesh elements to train the CNN. However, in this case one should be able to design a strategy to label elements automatically.}
%Again, notice that the performance could be improved by considering more data and using an architecture with more layers We also point out that our goal is not to optimize this process, but rather to show the importance of a classification step in the refinement procedure, and how CNNs can be employed for this purpose.
\section{Validation on a set of polygonal meshes}
In this section we compare the performance of the proposed algorithms. We consider four different coarse grids of the domain $(0,1)^2$: a grid of triangles, a Voronoi grid, a smoothed Voronoi grid obtained with  Polymesher \cite{talischi2012polymesher}, and a grid made of non-convex elements. In Figure~\ref{table: corse grids refined} these grids have been successively refined uniformly, i.e. each mesh element has been refined, for three times using the MP, the CNN-MP and the CNN-RP strategies. The final number of mesh elements is shown in Table \ref{table: mesh elements}. We observe that on average the MP strategy produced {\color{black}3 or 4 times} more elements than the CNN-RP strategy, and {\color{black}5 times} more than CNN-MP strategy.\\
\newline
\begin{table}
\centering
{\color{black}
\begin{tabular}{@{}c|cccc@{}}
\textbf{\# mesh elements} & triangles & Voronoi & smoothed Voronoi & non-convex \\ \hline
initial grid  & 32        & 9       & 10               & 14           \\
MP            & 6371      & 2719    & 3328             & 4502         \\
CNN-RP        & 2048      & 685     & 784              & 1574     \\
CNN-MP        & 1201      & 536     & 667              & 849          \\
\end{tabular}
}
\caption[Mesh elements of grids refined uniformly using different strategies.]
{Final number of elements for each mesh shown in Figure~\ref{table: corse grids refined}: a grid of triangles, a Voronoi grid, a smoothed Voronoi grid and a grid made of non-convex elements have been uniformly refined using the Midp-Point (MP), the CNN-enhanced Mid-Point (CNN-MP) and the CNN-enhanced Reference Polygon (CNN-RP) strategies. On average, the MP strategy produced {\color{black}3 or 4 times} more elements than the CNN-RP strategy, and {\color{black}5 times} more elements than CNN-MP strategy.}
\label{table: mesh elements}
\end{table}
\begin{table}[p!]
        \hspace{-2.5cm}
        \begin{tabular}{cM{35mm}M{35mm}M{35mm}M{35mm}}
             & initial grid & MP & CNN-RP & CNN-MP \\
             
             \rotatebox[origin=c]{90}{triangles}
             &
             \includegraphics[width = \linewidth]{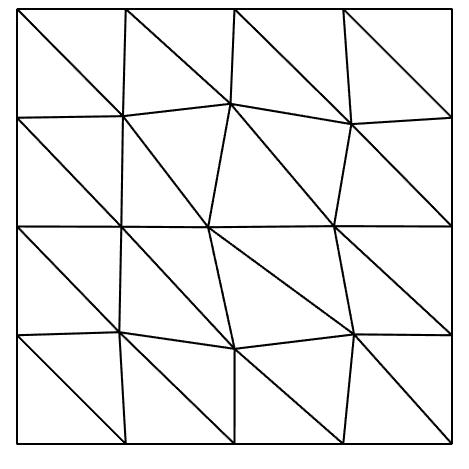} &
             \includegraphics[width=\linewidth]{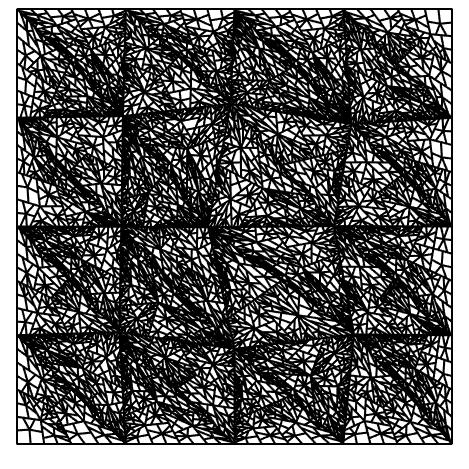}
             &
             \includegraphics[width=\linewidth]{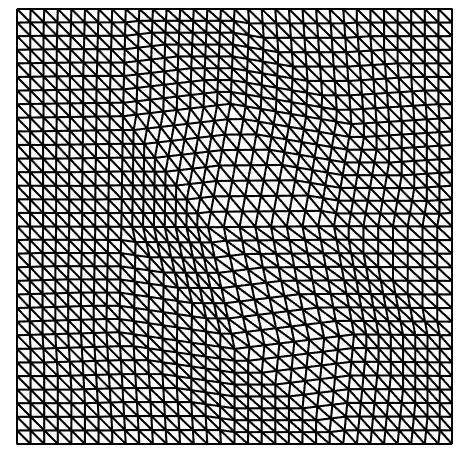}
             &
             \includegraphics[width=\linewidth]{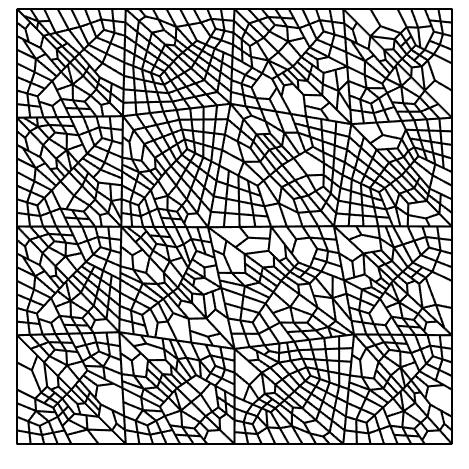}
             \\
            
            \rotatebox[origin=c]{90}{Voronoi}
            &
            \includegraphics[width = \linewidth]{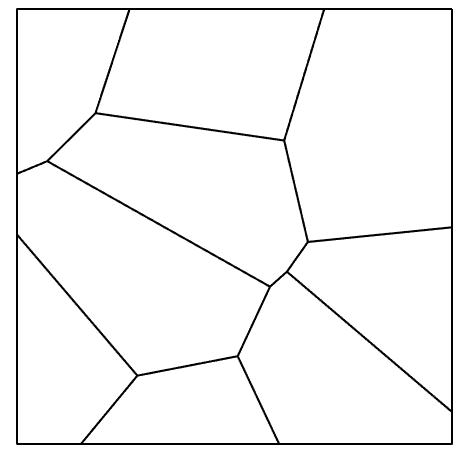}
            &
            \includegraphics[width=\linewidth]{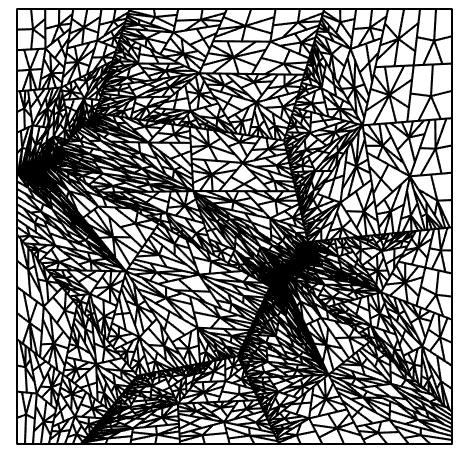}
            &
            \includegraphics[width=\linewidth]{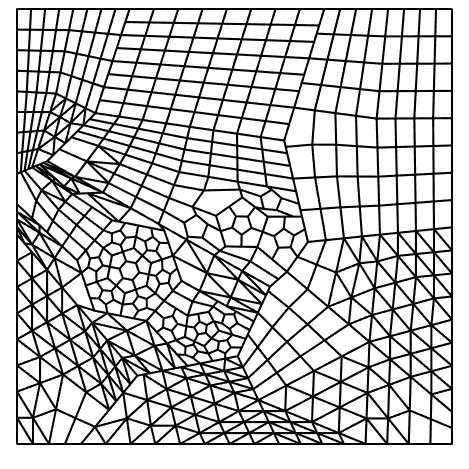}
            &
            \includegraphics[width=\linewidth]{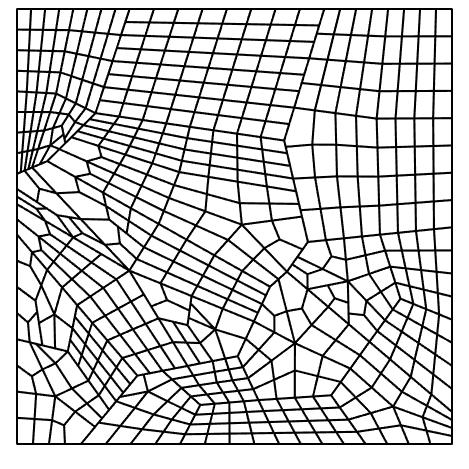}
            \\
            
            \rotatebox[origin=c]{90}{smoothed Voronoi}
            &
            \includegraphics[width = \linewidth]{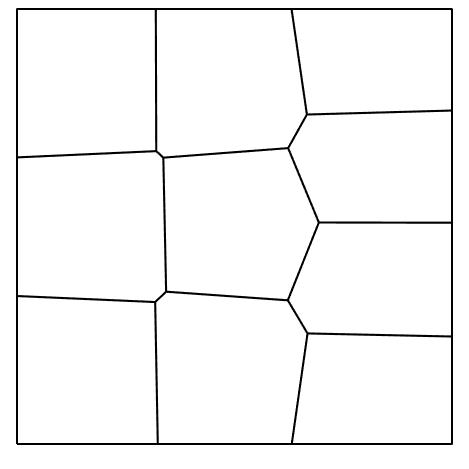}
            &
            \includegraphics[width=\linewidth]{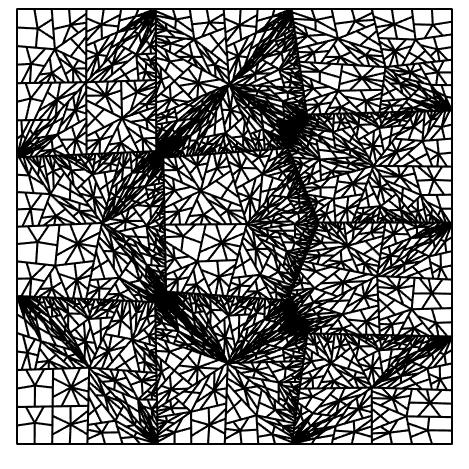}
            &
            \includegraphics[width=\linewidth]{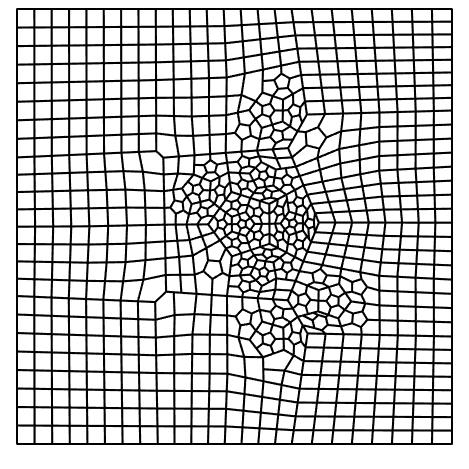}
            &
            \includegraphics[width=\linewidth]{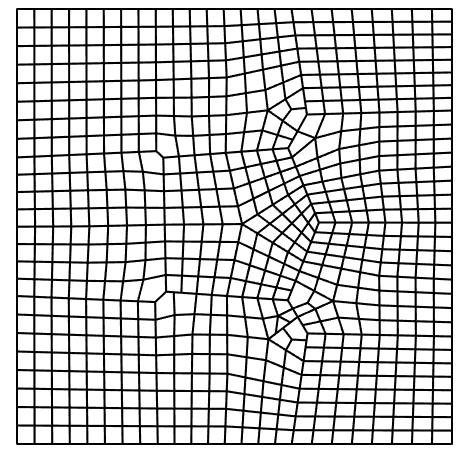}
            \\
            
            \rotatebox[origin=c]{90}{non-convex}
            &
            \includegraphics[width = \linewidth]{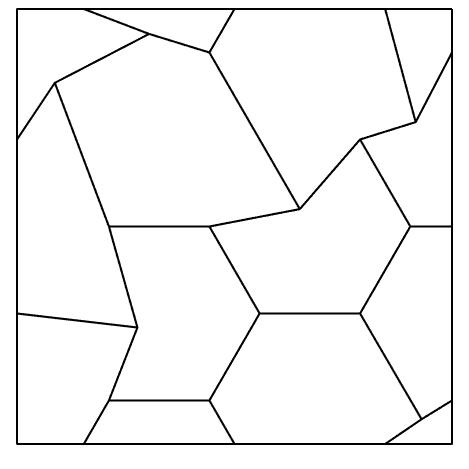}
            &
            \includegraphics[width=\linewidth]{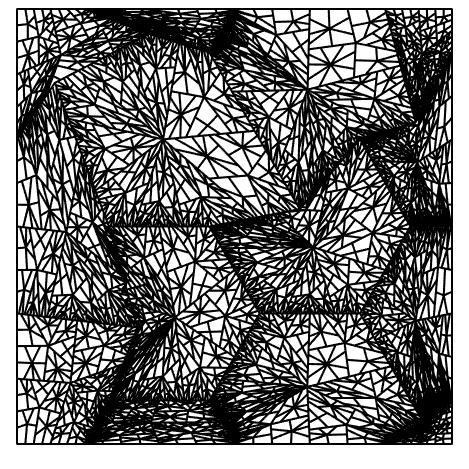}
            &
            \includegraphics[width=\linewidth]{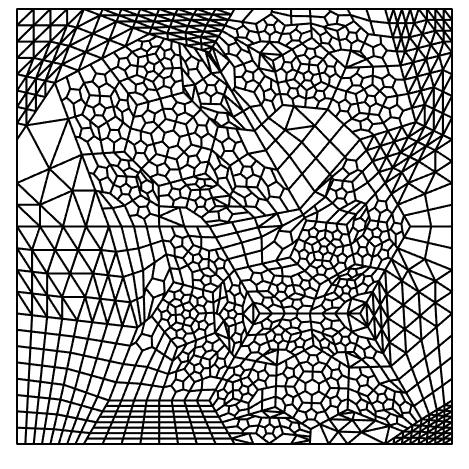}
            &
            \includegraphics[width=\linewidth]{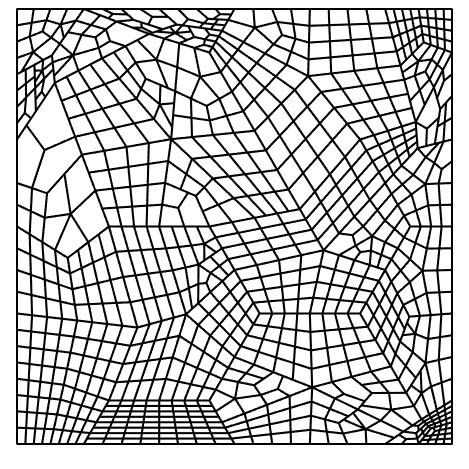}
            \\
        \end{tabular}
        \captionsetup{width=1.25\linewidth}
        \captionof{figure}[Grids refined uniformly using different strategies.]{In the first column, coarse grids of the domain $\Omega = (0,1)^2$: a grid of triangles, a Voronoi grid, a smoothed Voronoi grid, and a grid made of non-convex elements. Second to fourth columns: refined grids obtained after three steps of uniform refinement based on employing the MP (second column), the CNN-RP (third column) and the CNN-MP (fourth column) strategies. Each row corresponds to the same initial grid, while each column corresponds to the same refinement strategy.}
        \label{table: corse grids refined}
\end{table}
%%%%%%%%%%%%%%%%%%%%%%%%%%%%%%%%%%%%%%%%%%%%
In Figure~\ref{fig:all metrics} we show the computed quality metrics described in Section~3.3 on the grids of Figure~\ref{table: corse grids refined} (triangles, Voronoi, smoothed Voronoi, non-convex). Despite the fact that the performance are considerably grid dependent, the CNN-RP  strategy and the CNN-MP strategy seem to perform in a comparable way. Moreover, the CNN-RP and the CNN-MP strategies perform consistently better than the MP strategy, since their distributions are generally more concentrated toward the value 1.

%%%%%%%%%%%%%%%%%%%%%%%%%%%%%%%%%%%%%%%%%%%%%%%%%%%
%\begin{changemargin}{-3cm}{-3cm}

\begin{table}[p!]
        \vspace{-3.2cm}
        \hspace{-3cm}
        \begin{tabular}{cM{36mm}M{36mm}M{36mm}M{36mm}}
        & triangles & Voronoi & smoothed Voronoi & non-convex \\
        \rotatebox[origin=c]{90}{\small Uniformity Factor}&
        \includegraphics[page=1,width = \linewidth]{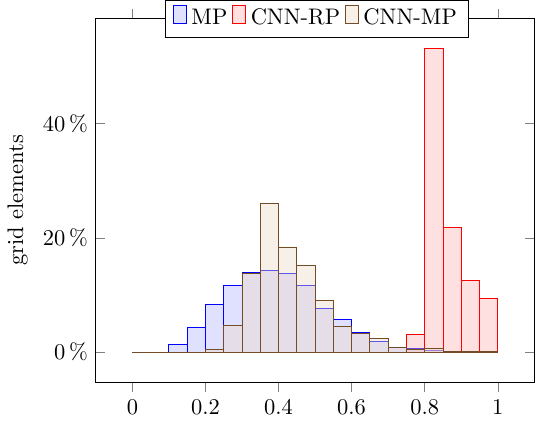}&
        \includegraphics[page=3,width = \linewidth]{images/quality_metrics.pdf}&
        \includegraphics[page=5,width = \linewidth]{images/quality_metrics.pdf}&
        \includegraphics[page=7,width = \linewidth]{images/quality_metrics.pdf}\\
        
        \rotatebox[origin=c]{90}{\small Circle Ratio}&
        \includegraphics[page=2,width = \linewidth]{images/quality_metrics.pdf}&
        \includegraphics[page=4,width = \linewidth]{images/quality_metrics.pdf}&
        \includegraphics[page=6,width = \linewidth]{images/quality_metrics.pdf}&
        \includegraphics[page=8,width = \linewidth]{images/quality_metrics.pdf}\\

        \rotatebox[origin=c]{90}{\small Area-Perimeter Ratio}&
        \includegraphics[page=9,width = \linewidth]{images/quality_metrics.pdf}&
        \includegraphics[page=11,width = \linewidth]{images/quality_metrics.pdf}&
        \includegraphics[page=13,width = \linewidth]{images/quality_metrics.pdf}&
        \includegraphics[page=15,width = \linewidth]{images/quality_metrics.pdf}\\

        \rotatebox[origin=c]{90}{\small Minimum Angle}&
        \includegraphics[page=10,width = \linewidth]{images/quality_metrics.pdf}&
        \includegraphics[page=12,width = \linewidth]{images/quality_metrics.pdf}&
        \includegraphics[page=14,width = \linewidth]{images/quality_metrics.pdf}&
        \includegraphics[page=16,width = \linewidth]{images/quality_metrics.pdf}\\

        \rotatebox[origin=c]{90}{\small Edge Ratio}&
        \includegraphics[page=17,width = \linewidth]{images/quality_metrics.pdf}&
        \includegraphics[page=19,width = \linewidth]{images/quality_metrics.pdf}&
        \includegraphics[page=21,width = \linewidth]{images/quality_metrics.pdf}&
        \includegraphics[page=23,width = \linewidth]{images/quality_metrics.pdf}\\

        \rotatebox[origin=c]{90}{\small Normalized Point Distance}&
        \includegraphics[page=18,width = \linewidth]{images/quality_metrics.pdf}&
        \includegraphics[page=20,width = \linewidth]{images/quality_metrics.pdf}&
        \includegraphics[page=22,width = \linewidth]{images/quality_metrics.pdf}&
        \includegraphics[page=24,width = \linewidth]{images/quality_metrics.pdf}\\
        
        \end{tabular}
        
                \captionsetup{width=1.25\linewidth}
        
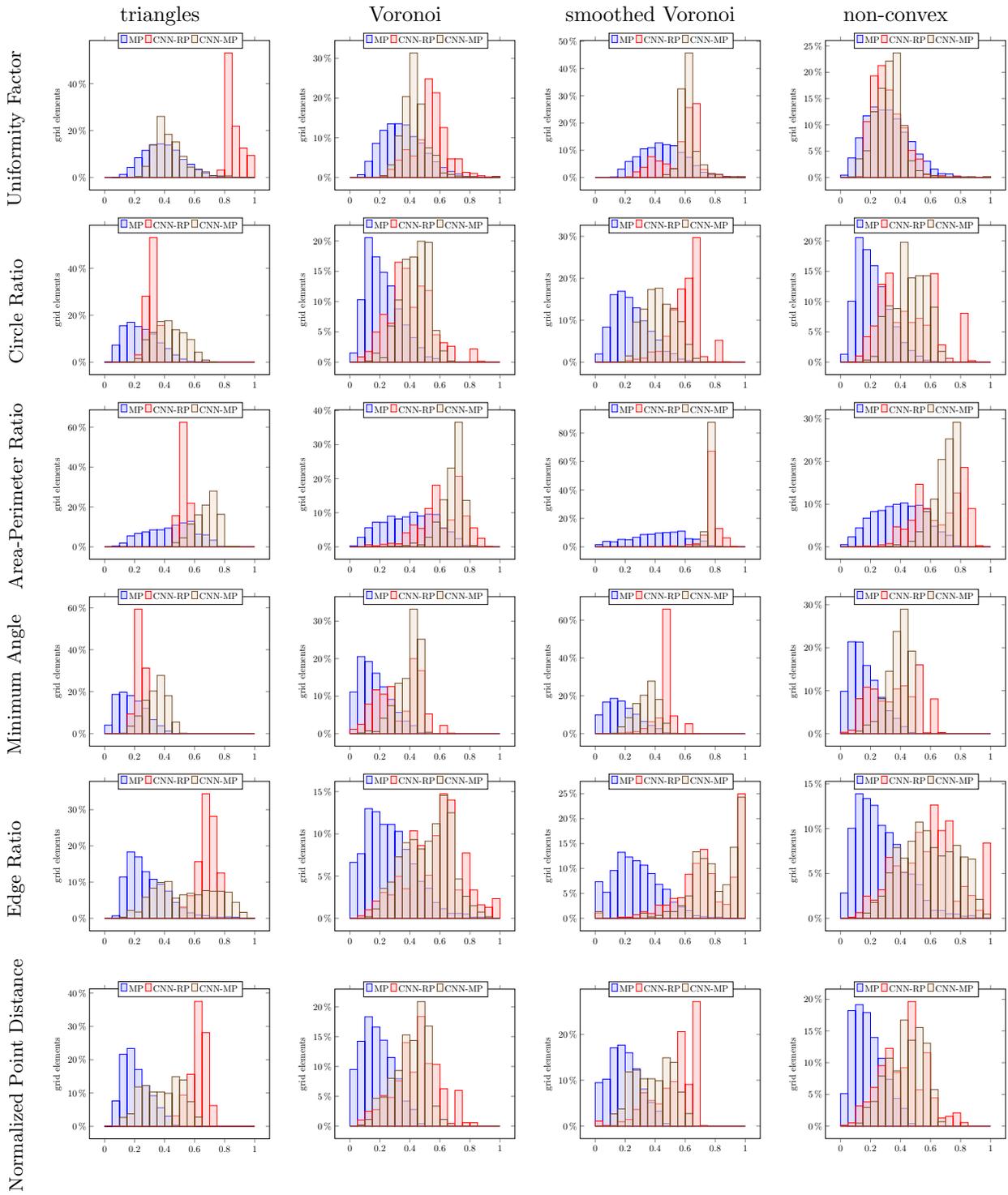
\captionof{figure}{Computed quality metrics (Uniformity Factor, Circle Ratio, Minimum Angle, Edge Ratio and Normalized Point Distance) for the refined grids reported in Figure~\ref{table: corse grids refined} (second to fourth column) and obtained based on employing different refinement strategies (MP, CNN-MP, CNN-RP).}
       \label{fig:all metrics}
\end{table}
\section{Testing CNN-based refinement strategies with PolyDG and Virtual Elements discretizations}
In this section we test the effectiveness of the proposed refinement strategies, to be used in combination with polygonal finite element discretizations. To this aim we consider PolyDG and Virtual Element discretizations of the following model problem: find $u \in H^1_0(\Omega)$ such that
\begin{equation}
 \label{eq:Poisson problem}
\int_\Omega \nabla u \cdot \nabla v = \int_\Omega f v \quad \forall v \in H^1_0(\Omega),
\end{equation}
with $f \in L^2(\Omega)$ a given forcing term.
The workflow is as follows:
\begin{enumerate}
    \item Generate a grid for $\Omega$.
    \item Compute numerically the solution of problem \eqref{eq:Poisson problem} using either the VEM \cite{beirao2013basic,beirao2014hitchhiker,beirao2016virtual,da2016mixed} or the PolyDG method  \cite{hesthaven2007nodal,bassi2012flexibility,antonietti2013hp,cangiani2014hp,antonietti2016review,cangiani2017hp}.
    \item Compute the error. In the VEM case the error is measured using the discrete $H^1_0$-like norm (see \cite{beirao2014hitchhiker}, for details), while in the PolyDG case the error is computed using the DG norm (see \cite{arnold2002unified,cockburn2012discontinuous}, for details)

    $$\|v\|_{\text{DG}}^{2}=\sum_{P}\|\nabla v\|_{L^{2}(P)}^{2}+ \sum_F \|\gamma^{1 / 2} \llbracket v \rrbracket \|_{L^{2}\left(F\right)}^{2}
    ,$$
    where $\gamma$ is the stabilization function (that depends on the discretization parameters and is chosen as in \cite{cangiani2014hp}), $P$ is a polygonal mesh element and $F$ is an element face. The jump operator $\llbracket \cdot \rrbracket$ is defined as in \cite{arnold2002unified}.
    
    \item

   Use the fixed fraction refinement strategy to refine a fraction $r$ of the number of elements. To refine the marked elements we employ one of the proposed strategies. Here, in order to investigate the effect of the proposed refinement strategies, we did not employ any a posteriori estimator of the error, but we computed element-wise the local error based on employing the exact solution.
\end{enumerate}

\subsection{Uniformly refined grids}
\begin{figure}[p!]
\vspace{-2cm}
\includegraphics[page=1,width = 0.49\linewidth]{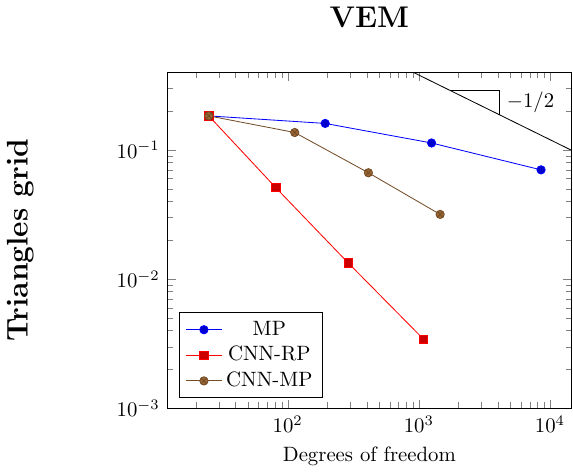}
\hfill
\includegraphics[page=2,width = 0.45\linewidth]{images/convergence.pdf}
\includegraphics[page=3,width = 0.49\linewidth]{images/convergence.pdf}
\hfill
\includegraphics[page=4,width = 0.45\linewidth]{images/convergence.pdf}
\includegraphics[page=5,width = 0.49\linewidth]{images/convergence.pdf}
\hfill
\includegraphics[page=6,width = 0.45\linewidth]{images/convergence.pdf}
\includegraphics[page=7,width = 0.49\linewidth]{images/convergence.pdf}
\hfill
\includegraphics[page=8,width = 0.45\linewidth]{images/convergence.pdf}
\caption{Test case of Section 6.1. Computed errors as a function of the number of degrees of freedom. Each row corresponds to the same initial grid (triangles, Voronoi, smoothed Voronoi, non-convex) refined uniformly with the proposed refinement strategies (MP, CNN-RP and CNN-MP), while each column corresponds to a different numerical method (VEM left and PolyDG right).}
\label{fig:VEM DG unif}
\end{figure}

When $r = 1$, the grid is refined uniformly, i.e. at each refinement step each mesh element is refined. The forcing term $f$ in \eqref{eq:Poisson problem} is selected in such a way that the exact solution is given by $u(x,y) = \sin(\pi x)\sin(\pi y).$
The grids obtained after three steps of uniform refinement are those already reported in Figure~\ref{table: corse grids refined}. In Figure~\ref{fig:VEM DG unif} we show the computed errors as a function of the number of degrees of freedom. We observe that the CNN-enhanced strategies (both MP and RP ones) outperform the plain MP rule. The difference is more evident for VEMs than for PolyDG approximations. This different sensitivity to mesh distortions may be due to the fact that VEM are hybrid methods with unknowns on the elements boundary.\\
{\color{black}Generating the binary image of a polygon has a computational cost that scales linearly with the number of edges of a polygon and with the squared root of the total number of pixels. Evaluating online the CNN has a computational cost that depends on the number pixels of the input image (and on the classification problem at hand). With our current implementation and for the considered benchmarks, classifying the shape of every mesh element takes on average approximately 3\% the time of solving the numerical problem over the refined grid. However, this ratio will decrease if meshes with more elements are considered, because except for particular cases the average number of edges of a mesh element will remain approximately constant and the number of pixels will remain constant, while solving the numerical problem has a cost which in general scales more than linearly with the number of degrees of freedom.}

\begin{figure}[p!]
\vspace{-2cm}
\includegraphics[page=25,width = 0.49\linewidth]{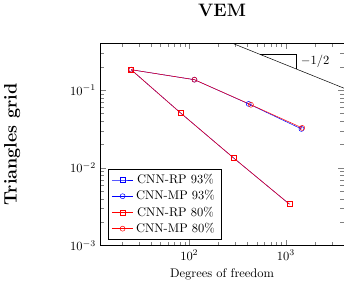}
\hfill
\includegraphics[page=26,width = 0.45\linewidth]{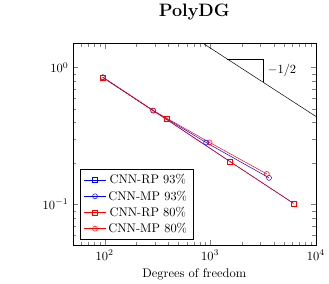}
\includegraphics[page=27,width = 0.49\linewidth]{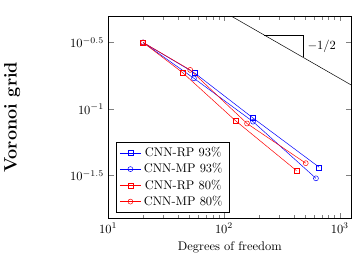}
\hfill
\includegraphics[page=28,width = 0.45\linewidth]{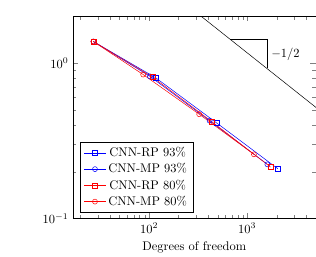}
\includegraphics[page=29,width = 0.49\linewidth]{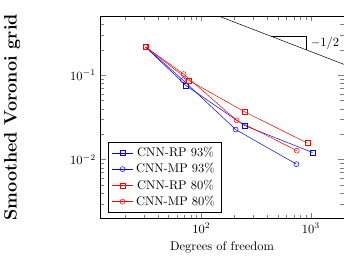}
\hfill
\includegraphics[page=30,width = 0.45\linewidth]{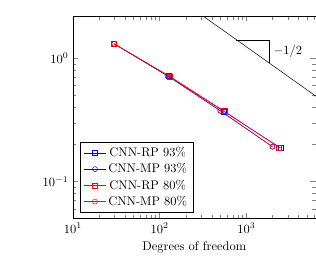}
\includegraphics[page=31,width = 0.49\linewidth]{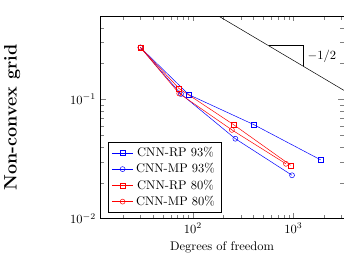}
\hfill
\includegraphics[page=32,width = 0.45\linewidth]{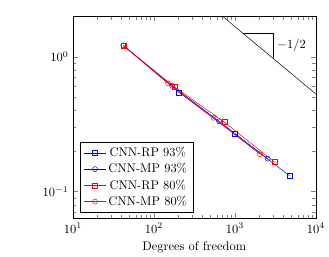}
\caption{{\color{black}Test case of Section 6.1. Performance comparison of the currently used CNN with accuracy 93\% and of a CNN with accuracy 80\%. Computed errors as a function of the number of degrees of freedom. Each row corresponds to the same initial grid (triangles, Voronoi, smoothed Voronoi, non-convex) uniformly refined with the CNN-RP and CNN-MP refinement strategies, respectively. Each column corresponds to a different numerical method (VEM left and PolyDG right).}}
\label{fig:net80-90}
\end{figure}
{\color{black} In Figure \ref{fig:net80-90} we compare the performance of the currently used CNN with accuracy 93\% and of a CNN with accuracy 80\%. As we can see, the PolyDG method does not seem very sensitive to the CNN accuracy, while the VEM seems more sensitive but performance seem not to always improve consistently over the selected grids. 

}
\subsection{Adaptively refined grids}
\begin{figure}[p!]
\vspace{-2cm}
\includegraphics[page=9,width = 0.49\linewidth]{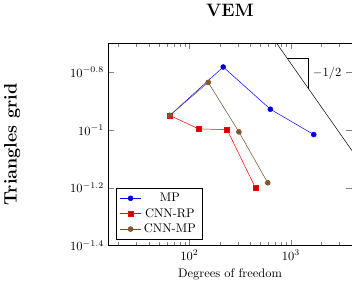}
\hfill
\includegraphics[page=10,width = 0.45\linewidth]{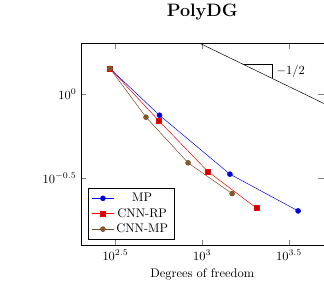}
\includegraphics[page=11,width = 0.49\linewidth]{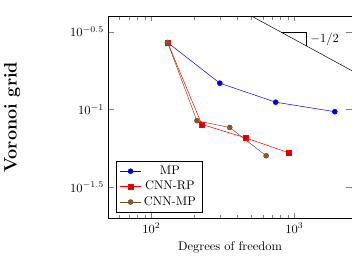}
\hfill
\includegraphics[page=12,width = 0.45\linewidth]{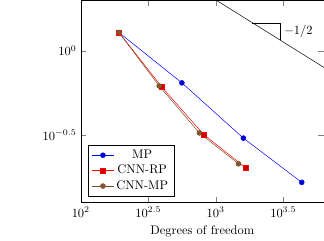}
\includegraphics[page=13,width = 0.49\linewidth]{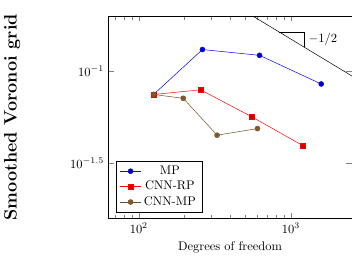}
\hfill
\includegraphics[page=14,width = 0.45\linewidth]{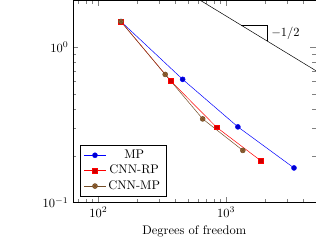}
\includegraphics[page=15,width = 0.49\linewidth]{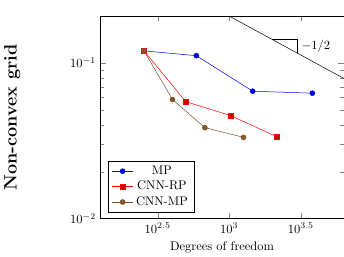}
\hfill
\includegraphics[page=16,width = 0.45\linewidth]{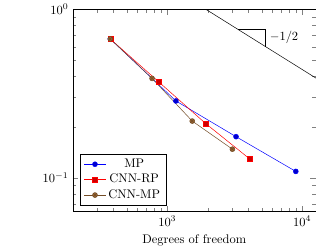}
\caption{Test case of Section 6.2. Computed errors as a function of the number of degrees of freedom. Each row corresponds to the same initial grid (triangles, Voronoi, smoothed Voronoi, non-convex) refined adaptively with a fixed fraction refinement criterion (refinement fraction $r$ set equal to 30\%) with different strategies (MP, CNN-RP and CNN-MP), while each column corresponds to a different numerical method (VEM left and PolyDG right).}
\label{fig:VEM DG adapt}
\end{figure}

In this case we selected $r = 0.3$. The forcing term $f$ in \eqref{eq:Poisson problem} is selected in such a way that the exact solution is $u(x,y) = (1 - e^{-10x})(x-1)\sin(\pi y)$, that exhibits a boundary layer along  $x = 0$. Figure~\ref{table: fine grids refined DG} shows the computed grids after three steps of refinement for the PolyDG case. Very similar grids have been obtained with Virtual Element discretizations.
    \begin{table}[p!]
        \hspace{-2.6cm}
        \begin{tabular}{cM{35mm}M{35mm}M{35mm}M{35mm}}
            & initial grid & MP & CNN-RP & CNN-MP \\
             
             \rotatebox[origin=c]{90}{triangles}
             &
             \includegraphics[width = \linewidth]{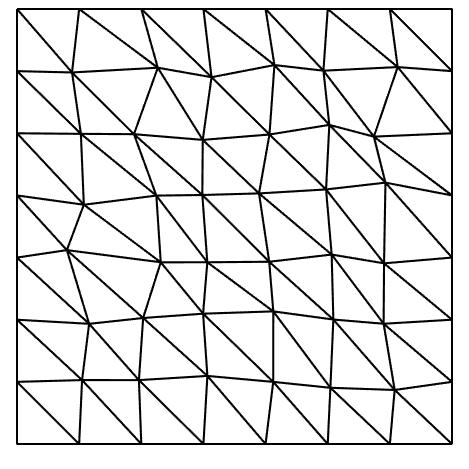} &
             \includegraphics[width=\linewidth]{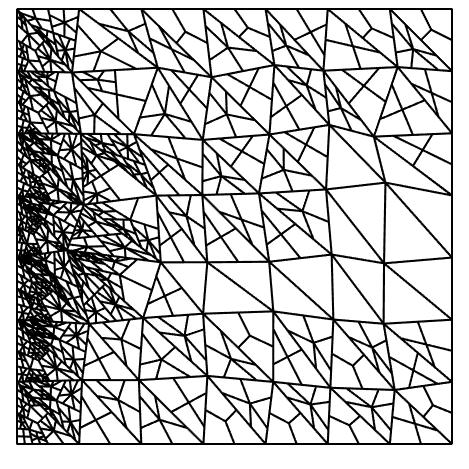}
             &
             \includegraphics[width=\linewidth]{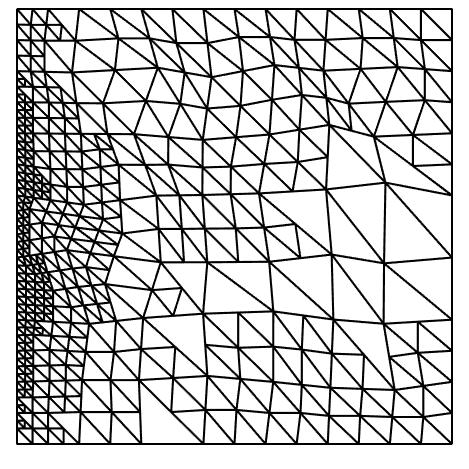}
             &
             \includegraphics[width=\linewidth]{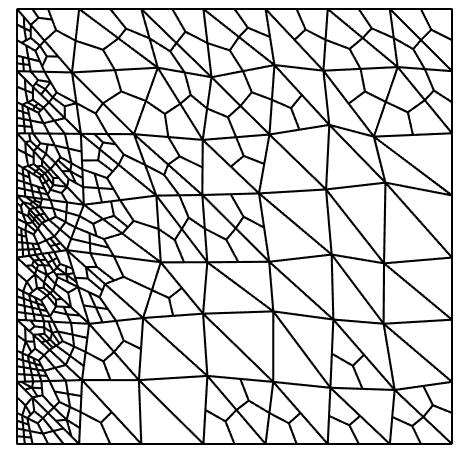}
             \\
            
            \rotatebox[origin=c]{90}{Voronoi}
            &
            \includegraphics[width = \linewidth]{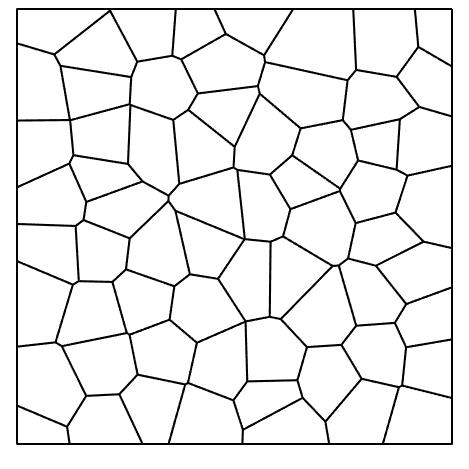}
            &
            \includegraphics[width=\linewidth]{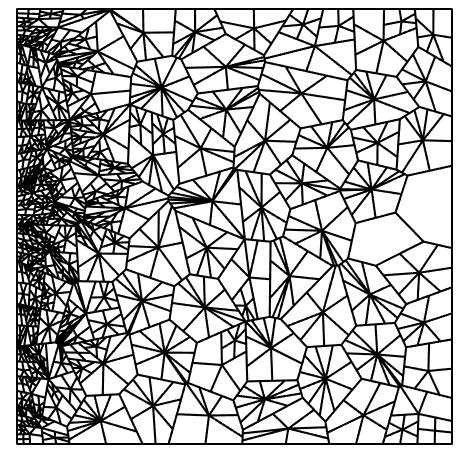}
            &
            \includegraphics[width=\linewidth]{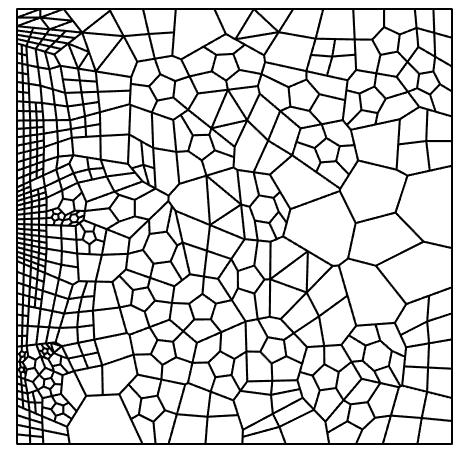}
            &
            \includegraphics[width=\linewidth]{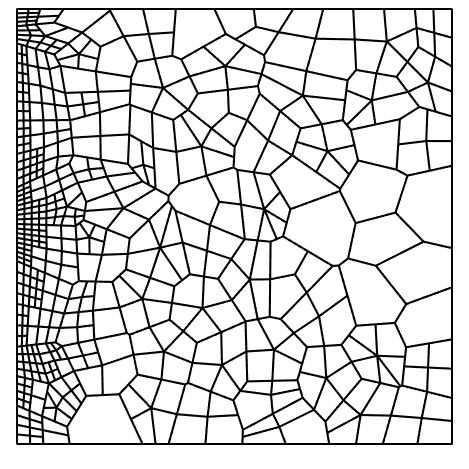}
            \\
            
            \rotatebox[origin=c]{90}{smoothed Voronoi}
            &
            \includegraphics[width = \linewidth]{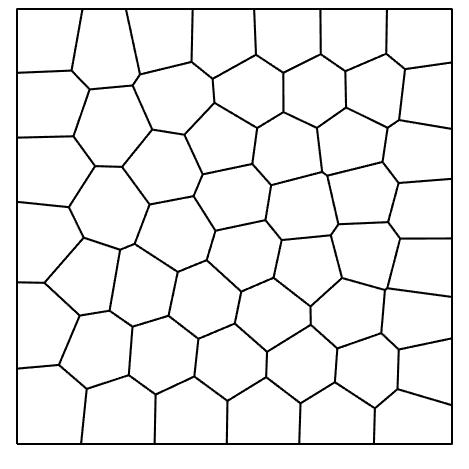}
            &
            \includegraphics[width=\linewidth]{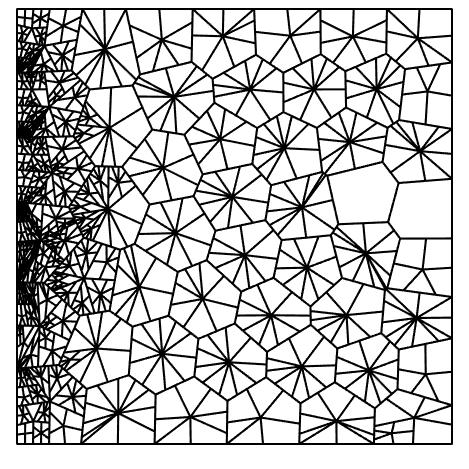}
            &
            \includegraphics[width=\linewidth]{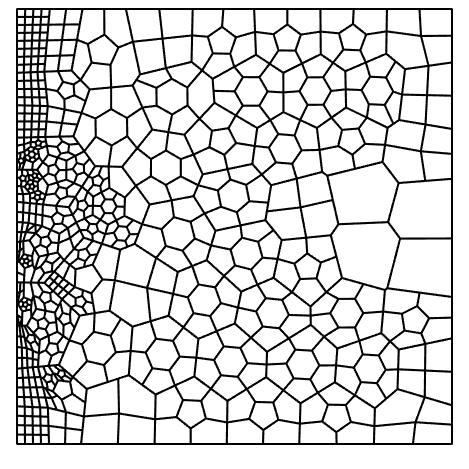}
            &
            \includegraphics[width=\linewidth]{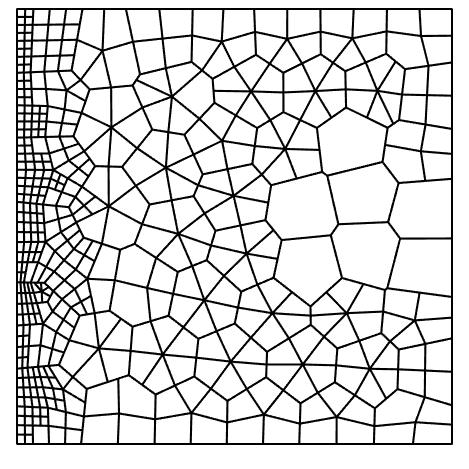}
            \\
            
            \rotatebox[origin=c]{90}{non-convex}
            &
            \includegraphics[width = \linewidth]{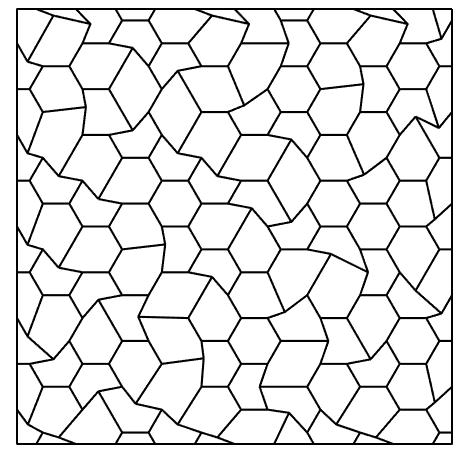}
            &
            \includegraphics[width=\linewidth]{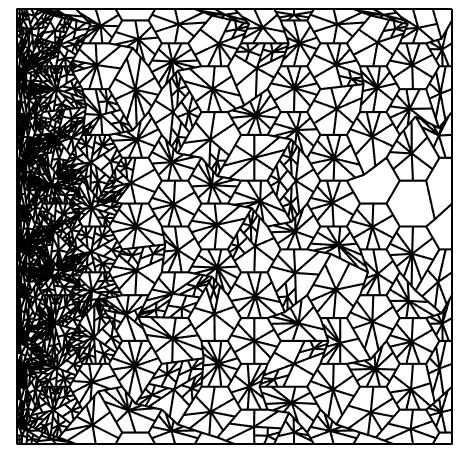}
            &
            \includegraphics[width=\linewidth]{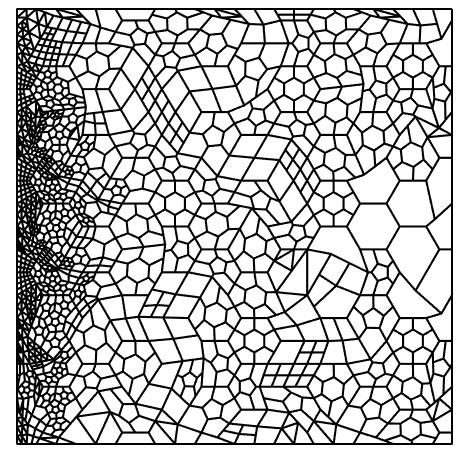}
            &
            \includegraphics[width=\linewidth]{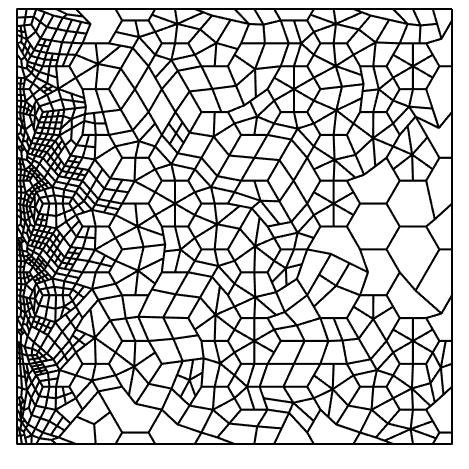}
            \\
        \end{tabular}
        \captionsetup{width=1.25\linewidth}
        \captionof{figure}[Grids refined adaptively using the DG method.]{Adaptively refined grids for the test case of Section 6.2. Each row corresponds to the same initial grid (triangles, Voronoi, smoothed Voronoi, non-convex), while the second-fourth columns correspond to the different refinement strategies (MP, CNN-RP, CNN-MP). Three successively adaptive refinement steps have been performed, with a fixed fraction refinement criterion (refinement fraction $r$ set equal to 30\%).}
        \label{table: fine grids refined DG}
    \end{table}
\\
In Figure~\ref{fig:VEM DG adapt} we show the computed errors as a function of the number of degrees of freedom for both Virtual Element and PolyDG discretizations. 
The CNN-enhanced strategies (both MP and RP ones) outperform the plain MP rule. The difference is more evident for VEMs than for PolyDG approximations.

{\color{black}
\subsection{Application to an advection-diffusion problem}
We now consider the following advection-diffusion problem:
$$
\left\{\begin{array}{ll} - \Delta u + \operatorname{div} ( \boldsymbol{\beta} u)
=f & \text { in } \Omega \\
u=0 & \text { on } \partial \Omega
\end{array}\right.
$$
where $\Omega = (0,1)^2$ and $\boldsymbol{\beta} = [1\ 0]^T$ is a constant velocity filed. We choose the forcing term $\mathbf{f}$ in such a way that the exact solution is $(x-1)(1-\exp(2x))\sin(\pi y)$,
%and solve it using PolyDG in the uniform refinement case over the grids of Figure \ref{table: corse grids refined}. 
and solve this problem using PolyDG method, both in the uniform refinement case (coarse grids in the first column of Figure \ref{table: corse grids refined}), and in the adaptive refinement case with fixed refinement fraction of 30\% (fine grids in the first column of Figure \ref{table: fine grids refined DG}).
\begin{table}[p!]
        \vspace{-2cm}
        % \hspace{-2.2cm}
        \begin{tabular}{cM{45mm}M{45mm}}
             & uniform refinement & adaptive refinement \\
             
             \rotatebox[origin=c]{90}{triangles}
             &
             \includegraphics[page=33,width = \linewidth]{images/convergence.pdf} &
             \includegraphics[page=34,width = \linewidth]{images/convergence.pdf}
             \\
            
            \rotatebox[origin=c]{90}{Voronoi}
            &
            \includegraphics[page=35,width = \linewidth]{images/convergence.pdf}
            &
            \includegraphics[page=36,width = \linewidth]{images/convergence.pdf}
            \\
            
            \rotatebox[origin=c]{90}{smoothed Voronoi}
            &
            \includegraphics[page=37,width = \linewidth]{images/convergence.pdf}
            &
            \includegraphics[page=38,width = \linewidth]{images/convergence.pdf}
            \\
            
            \rotatebox[origin=c]{90}{non-convex}
            &
            \includegraphics[page=39,width = \linewidth]{images/convergence.pdf}
            &
           \includegraphics[page=40,width = \linewidth]{images/convergence.pdf}
            \\
        \end{tabular}
        \captionsetup{width=1.25\linewidth}
        
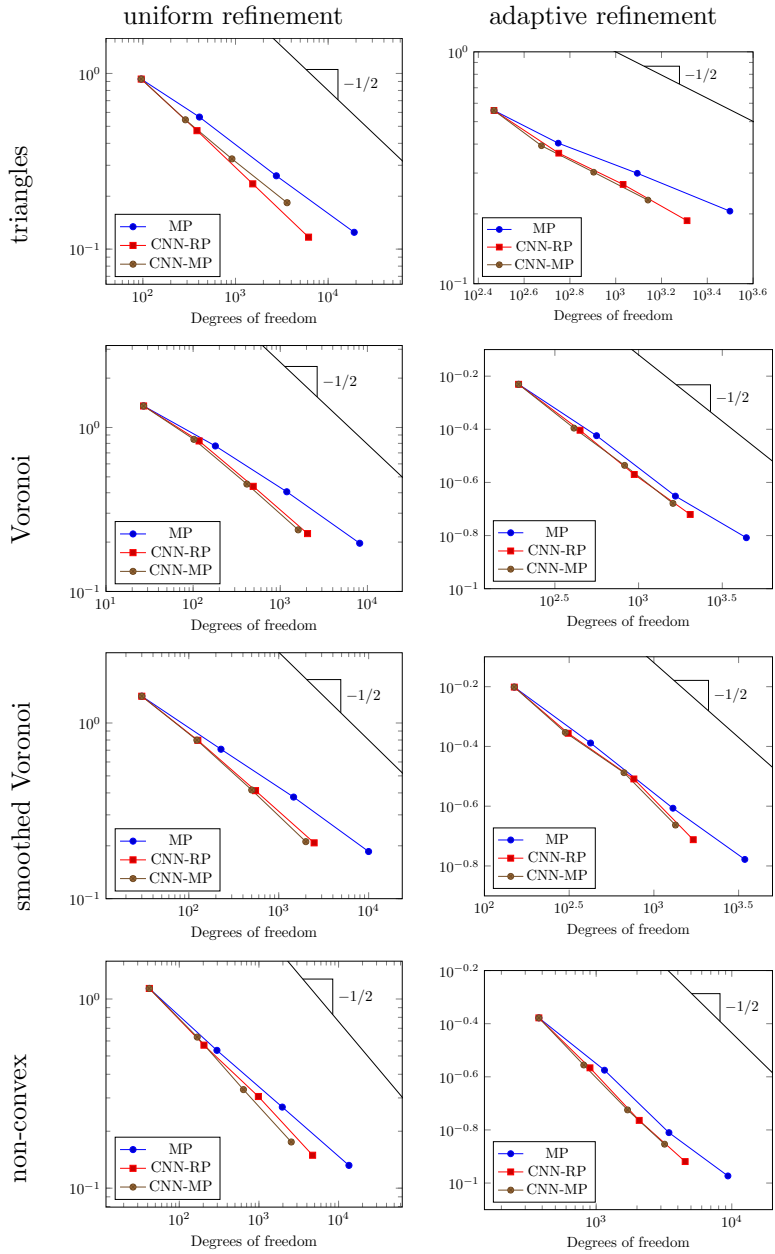
\captionof{figure}{\color{black}Advection-diffusion problem of Section 6.3. Computed errors using the DG norm as a function of the number of degrees of  freedom.   Each  row  corresponds  to  the  same  initial  grid  (triangles,  Voronoi,  smoothed Voronoi, non-convex) refined uniformly with the proposed refinement strategies (MP, CNN-RP and CNN-MP). First column: uniform refinement case over the coarse grids in the first column of Figure \ref{table: corse grids refined}. Second column: adaptive refinement case with fixed refinement fraction of 30\% of mesh elements over the fine grids in the first column of Figure \ref{table: fine grids refined DG}.}
        \label{fig:avection}
\end{table}
In Figure \ref{fig:avection} we show the computed errors using the DG norm as a function of the number of degrees of freedom. As we can see, the CNN-enhanced strategies (both MP and RP ones) perform better than the plain MP rule.}

{\color{black}
\subsection{Application to the Stokes problem}
We now consider the following Stokes problem, which describes the motion of an incompressible viscous flow:
$$
\left\{\begin{array}{ll} \partial_t \mathbf{u}
- \Delta \mathbf{u}+\nabla p=\mathbf{f} & \text { in } \Omega \\
\operatorname{div} \mathbf{u}=0 & \text { in } \Omega \\
\mathbf{u}=\mathbf{0} & \text { on } \partial \Omega
\end{array}\right.
$$
where $\Omega = (0,1)^2$, $\mathbf{u}: \Omega \rightarrow  \mathbb{R}^2$ is the velocity, $p: \Omega \rightarrow  \mathbb{R}$ is the pressure, $\mathbf{f}: \Omega \rightarrow  \mathbb{R}^2$ is the forcing term. We choose $\mathbf{f}$ in such a way that the exact solution is
$$
\boldsymbol{u}=\left[\begin{array}{l}
-\cos (2 \pi x) \sin (2 \pi y) \\
\sin (2 \pi x) \cos (2 \pi y)
\end{array}\right], \quad p=1-e^{-x(x-1)(x-0.5)^{2}-y(y-1)(y-0.5)^{2}},$$
%and solve it using PolyDG in the uniform refinement case over the grids of Figure \ref{table: corse grids refined}. 
and solve it using PolyDG method, both in the uniform refinement (coarse grids in the first column of Figure \ref{table: corse grids refined}), and in the adaptive refinement case with fixed refinement fraction of 40\% (fine grids in the first column of Figure \ref{table: fine grids refined DG}).
\begin{table}[p!]
        \vspace{-2cm}
        % \hspace{-2.2cm}
        \begin{tabular}{cM{45mm}M{45mm}}
             & uniform refinement & adaptive refinement \\
             
             \rotatebox[origin=c]{90}{triangles}
             &
             \includegraphics[page=17,width = \linewidth]{images/convergence.pdf} &
             \includegraphics[page=18,width = \linewidth]{images/convergence.pdf}
             \\
            
            \rotatebox[origin=c]{90}{Voronoi}
            &
            \includegraphics[page=19,width = \linewidth]{images/convergence.pdf}
            &
            \includegraphics[page=20,width = \linewidth]{images/convergence.pdf}
            \\
            
            \rotatebox[origin=c]{90}{smoothed Voronoi}
            &
            \includegraphics[page=21,width = \linewidth]{images/convergence.pdf}
            &
            \includegraphics[page=22,width = \linewidth]{images/convergence.pdf}
            \\
            
            \rotatebox[origin=c]{90}{non-convex}
            &
            \includegraphics[page=23,width = \linewidth]{images/convergence.pdf}
            &
           \includegraphics[page=24,width = \linewidth]{images/convergence.pdf}
            \\
        \end{tabular}
        \captionsetup{width=1.25\linewidth}
        
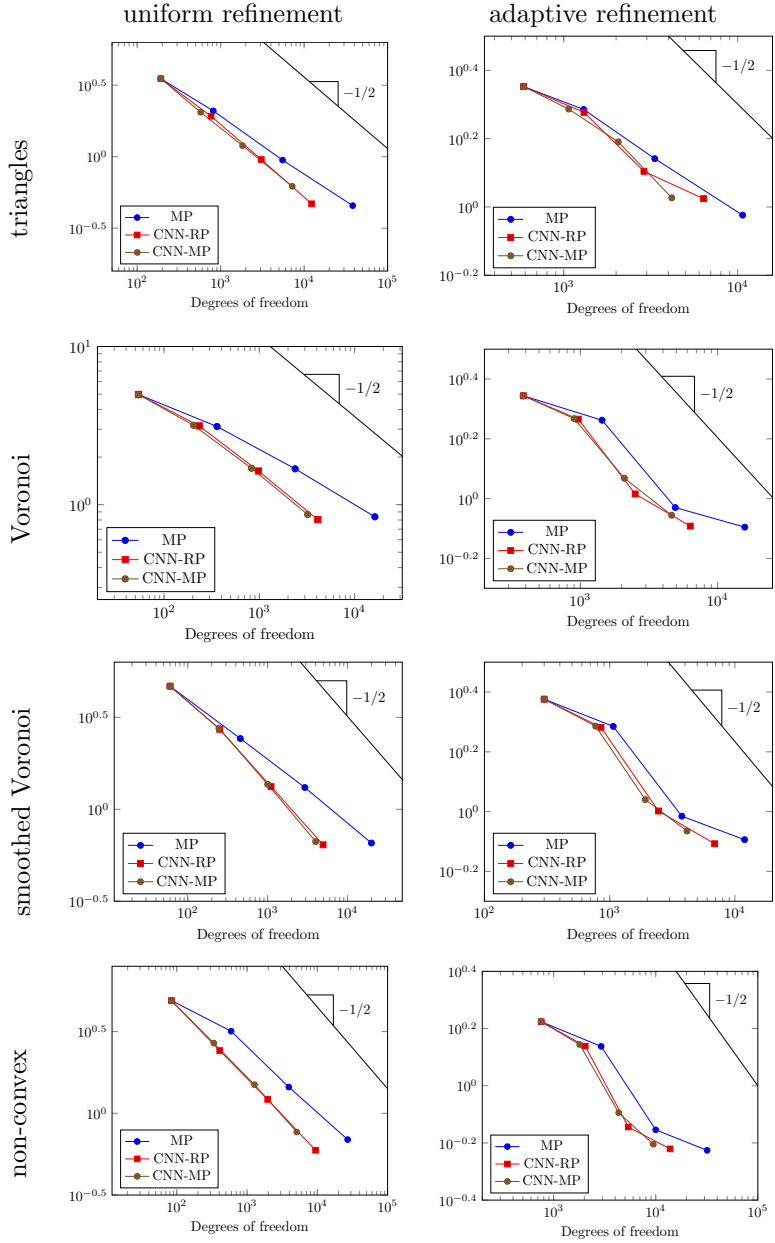
\captionof{figure}{\color{black}Stokes problem of Section 6.4. Computed errors in the $H^1$-broken norm of the velocity and $L^2$-norm of the pressure as a function of the number of degrees of  freedom.   Each  row  corresponds  to  the  same  initial  grid  (triangles,  Voronoi,  smoothed Voronoi, non-convex) refined uniformly with the proposed refinement strategies (MP, CNN-RP and CNN-MP). First column: uniform refinement case over the coarse grids in the first column of Figure \ref{table: corse grids refined}. Second column: adaptive refinement case with fixed refinement fraction of 40\% of mesh elements over the fine grids in the first column of Figure \ref{table: fine grids refined DG}.}
        \label{fig:Stokes unif adapt}
\end{table}
In Figure \ref{fig:Stokes unif adapt} we show the computed errors in the $H^1$-broken norm of the velocity and $L^2$-norm of the pressure as a function of the number of degrees of freedom. As it is clear from the results of Figure \ref{fig:Stokes unif adapt}, the CNN-enhanced strategies (both MP and RP ones) perform better than the plain MP rule.}

\section{Conclusions}
In this work, we propose a new paradigm based on CNNs to enhance both existing refinement criteria and new refinement procedures, withing polygonal finite element discretizations of partial differential equations. In particular, we introduced two refinement strategies for polygonal elements, named ``CNN-RP strategy" and ``CNN-MP strategy". The former proposes ad-hoc refinement strategies based on reference polygons, while the latter is an improved version of the known MP strategy. These strategies exploit a CNN to suitably classify the ``shape" of a polygon in order to later apply an ad-hoc refinement strategy. This approach has the advantage to be made of interchangeable pieces: any algorithm can be employed to classify mesh elements, as well as any refinement strategy can be employed to refine a polygon with a given label.\\
We have shown that correctly classifying elements' shape based on employing CNNs can improve consistently and significantly the quality of the grids and the accuracy of polygonal finite element methods employed for the discretization. Specifically, this has been measured in terms of less elements produced on average at each refinement step, in terms of improved quality of the mesh elements according to different quality metrics, and in terms of improved accuracy using numerical methods such as PolyDG methods and VEMs. These results show that classifying correctly the shape of a polygonal element plays a key role in which refinement strategy to choose, allowing to extend and to boost existing strategies. Moreover, this classification step has a very limited computational cost when using a pre-trained CNN. The latter can be made off-line once and for all, independently of the differential model under consideration.\\
As noticed in \cite{droniou2020interplay}, having skewed meshes can sometimes be beneficial to polytopal methods with hybrid unknowns, if the skewness is compatible with the diffusion anisotropy. We believe that our method may be extended so as to generate distorted elements by including the anisotropy directions into the shape recognition phase. Assume to dispose of a refinement strategy for ``squares", which consists in diving the element into other 4 ``squares", and for ``rectangles", which consists in diving the element into 2 ``squares". Assume you have to refine square $(0,1)^2$ and that the anisotropy directions are aligned with the axes. Before processing the shape of the element, you may skew the element along the anisotropy directions according to the magnitude of anisotropy, therefore turning the square into a rectangle. In this way the original square will be classified as ``rectangle" and divided into 2 rectangles skewed in the desired direction. When applying the same strategy on the new rectangles, they will be deformed into ``squares", and therefore refined into 4 elements skewed in the desired direction.\\
In terms of future research lines, we plan to extend these algorithms to three dimensional polyhedral grids. The CNN architecture is naturally designed to handle three dimensional images, while the design of effective refinement strategies in three dimensions is under investigation.
%%%%

\section{Acknowledgements}
Funding: P. F. Antonietti has been partially supported by the Ministero dell’Istruzione, dell’Università e della Ricerca [PRIN grant number 201744KLJL]. P. F. Antonietti is member of INDAM-GNCS.

\clearpage
\bibliography{main_bibfile}

\begin{thebibliography}{10}
\expandafter\ifx\csname url\endcsname\relax
  \def\url#1{\texttt{#1}}\fi
\expandafter\ifx\csname urlprefix\endcsname\relax\def\urlprefix{URL }\fi
\expandafter\ifx\csname href\endcsname\relax
  \def\href#1#2{#2} \def\path#1{#1}\fi

\bibitem{hyman1997numerical}
J.~Hyman, M.~Shashkov, S.~Steinberg, The numerical solution of diffusion
  problems in strongly heterogeneous non-isotropic materials, Journal of
  Computational Physics 132~(1) (1997) 130--148.

\bibitem{brezzi2005family}
F.~Brezzi, K.~Lipnikov, V.~Simoncini, A family of mimetic finite difference
  methods on polygonal and polyhedral meshes, Mathematical Models and Methods
  in Applied Sciences 15~(10) (2005) 1533--1551.

\bibitem{brezzi2005convergence}
F.~Brezzi, K.~Lipnikov, M.~Shashkov, Convergence of the mimetic finite
  difference method for diffusion problems on polyhedral meshes, SIAM Journal
  on Numerical Analysis 43~(5) (2005) 1872--1896.

\bibitem{da2014mimetic}
L.~Beirao~da Veiga, K.~Lipnikov, G.~Manzini, The mimetic finite difference
  method for elliptic problems, Vol.~11, Springer, 2014.

\bibitem{cockburn2008superconvergent}
B.~Cockburn, B.~Dong, J.~Guzm{\'a}n, A superconvergent ldg-hybridizable
  galerkin method for second-order elliptic problems, Mathematics of
  Computation 77~(264) (2008) 1887--1916.

\bibitem{cockburn2009superconvergent}
B.~Cockburn, J.~Guzm{\'a}n, H.~Wang, Superconvergent discontinuous galerkin
  methods for second-order elliptic problems, Mathematics of Computation
  78~(265) (2009) 1--24.

\bibitem{cockburn2009unified}
B.~Cockburn, J.~Gopalakrishnan, R.~Lazarov, Unified hybridization of
  discontinuous galerkin, mixed, and continuous galerkin methods for second
  order elliptic problems, SIAM Journal on Numerical Analysis 47~(2) (2009)
  1319--1365.

\bibitem{cockburn2010projection}
B.~Cockburn, J.~Gopalakrishnan, F.-J. Sayas, A projection-based error analysis
  of hdg methods, Mathematics of Computation 79~(271) (2010) 1351--1367.

\bibitem{hesthaven2007nodal}
J.~S. Hesthaven, T.~Warburton, Nodal discontinuous Galerkin methods:
  algorithms, analysis, and applications, Springer Science \& Business Media,
  2007.

\bibitem{bassi2012flexibility}
F.~Bassi, L.~Botti, A.~Colombo, D.~A. Di~Pietro, P.~Tesini, On the flexibility
  of agglomeration based physical space discontinuous galerkin discretizations,
  Journal of Computational Physics 231~(1) (2012) 45--65.

\bibitem{antonietti2013hp}
P.~F. Antonietti, S.~Giani, P.~Houston, hp-version composite discontinuous
  galerkin methods for elliptic problems on complicated domains, SIAM Journal
  on Scientific Computing 35~(3) (2013) A1417--A1439.

\bibitem{cangiani2014hp}
A.~Cangiani, E.~H. Georgoulis, P.~Houston, hp-version discontinuous galerkin
  methods on polygonal and polyhedral meshes, Mathematical Models and Methods
  in Applied Sciences 24~(10) (2014) 2009--2041.

\bibitem{antonietti2016review}
P.~F. Antonietti, A.~Cangiani, J.~Collis, Z.~Dong, E.~H. Georgoulis, S.~Giani,
  P.~Houston, Review of discontinuous galerkin finite element methods for
  partial differential equations on complicated domains, in: Building bridges:
  connections and challenges in modern approaches to numerical partial
  differential equations, Springer, 2016, pp. 281--310.

\bibitem{cangiani2017hp}
A.~Cangiani, Z.~Dong, E.~H. Georgoulis, P.~Houston, hp-Version discontinuous
  Galerkin methods on polygonal and polyhedral meshes, Springer, 2017.

\bibitem{beirao2013basic}
L.~Beir{\~a}o~da Veiga, F.~Brezzi, A.~Cangiani, L.~D. Manzini,
  GianM.and~Marini, A.~Russo, Basic principles of virtual element methods,
  Mathematical Models and Methods in Applied Sciences 23~(01) (2013) 199--214.

\bibitem{beirao2014hitchhiker}
L.~Beir{\~a}o~da Veiga, F.~Brezzi, L.~D. Marini, A.~Russo, The hitchhiker's
  guide to the virtual element method, Mathematical models and methods in
  applied sciences 24~(08) (2014) 1541--1573.

\bibitem{beirao2016virtual}
L.~Beir{\~a}o~da Veiga, F.~Brezzi, L.~Marini, A.~Russo, Virtual element method
  for general second-order elliptic problems on polygonal meshes, Mathematical
  Models and Methods in Applied Sciences 26~(04) (2016) 729--750.

\bibitem{da2016mixed}
L.~Beirao~da Veiga, F.~Brezzi, L.~D. Marini, A.~Russo, Mixed virtual element
  methods for general second order elliptic problems on polygonal meshes,
  ESAIM: Mathematical Modelling and Numerical Analysis 50~(3) (2016) 727--747.

\bibitem{di2014arbitrary}
D.~A. Di~Pietro, A.~Ern, S.~Lemaire, An arbitrary-order and compact-stencil
  discretization of diffusion on general meshes based on local reconstruction
  operators, Computational Methods in Applied Mathematics 14~(4) (2014)
  461--472.

\bibitem{di2015hybrid}
D.~A. Di~Pietro, A.~Ern, A hybrid high-order locking-free method for linear
  elasticity on general meshes, Computer Methods in Applied Mechanics and
  Engineering 283 (2015) 1--21.

\bibitem{di2015hybrid2}
D.~A. Di~Pietro, A.~Ern, Hybrid high-order methods for variable-diffusion
  problems on general meshes, Comptes Rendus Math{\'e}matique 353~(1) (2015)
  31--34.

\bibitem{di2016review}
D.~A. Di~Pietro, A.~Ern, S.~Lemaire, A review of hybrid high-order methods:
  formulations, computational aspects, comparison with other methods, in:
  Building bridges: connections and challenges in modern approaches to
  numerical partial differential equations, Springer, 2016, pp. 205--236.

\bibitem{di2019hybrid}
D.~A. Di~Pietro, J.~Droniou, The Hybrid High-Order method for polytopal meshes,
  Vol.~19, Springer, 2019.

\bibitem{attene2019benchmark}
M.~Attene, S.~Biasotti, S.~Bertoluzza, D.~Cabiddu, M.~Livesu, G.~Patan{\`e},
  M.~Pennacchio, D.~Prada, M.~Spagnuolo, Benchmark of polygon quality metrics
  for polytopal element methods, arXiv preprint arXiv:1906.01627.

\bibitem{di2021polyhedral}
D.~A. Di~Pietro, L.~Formaggia, R.~Masson, et~al., Polyhedral methods in
  geosciences.

\bibitem{burman2021convergence}
E.~Burman, O.~Duran, A.~Ern, M.~Steins, Convergence analysis of hybrid
  high-order methods for the wave equation, Journal of Scientific Computing
  87~(3) (2021) 1--30.

\bibitem{lai2016recursive}
M.-J. Lai, G.~Slavov, On recursive refinement of convex polygons, Computer
  Aided Geometric Design 45 (2016) 83--90.

\bibitem{hoshina2018simple}
T.~Y. Hoshina, I.~F. Menezes, A.~Pereira, A simple adaptive mesh refinement
  scheme for topology optimization using polygonal meshes, Journal of the
  Brazilian Society of Mechanical Sciences and Engineering 40~(7) (2018) 348.

\bibitem{berrone2021refinement}
S.~Berrone, A.~Borio, A.~D'Auria, Refinement strategies for polygonal meshes
  applied to adaptive vem discretization, Finite Elements in Analysis and
  Design 186 (2021) 103502.

\bibitem{chan1998agglomeration}
T.~F. Chan, J.~Xu, L.~Zikatanov, An agglomeration multigrid method for
  unstructured grids, Contemporary Mathematics 218 (1998) 67--81.

\bibitem{antonietti2020agglomeration}
P.~F. Antonietti, P.~Houston, G.~Pennesi, E.~S{\"u}li, An agglomeration-based
  massively parallel non-overlapping additive schwarz preconditioner for
  high-order discontinuous galerkin methods on polytopic grids, Mathematics of
  Computation.

\bibitem{lecun2015deep}
Y.~LeCun, Y.~Bengio, G.~Hinton, Deep learning, nature 521~(7553) (2015)
  436--444.

\bibitem{raissi2019physics}
M.~Raissi, P.~Perdikaris, G.~E. Karniadakis, Physics-informed neural networks:
  A deep learning framework for solving forward and inverse problems involving
  nonlinear partial differential equations, Journal of Computational Physics
  378 (2019) 686--707.

\bibitem{raissi2018hidden}
M.~Raissi, G.~E. Karniadakis, Hidden physics models: Machine learning of
  nonlinear partial differential equations, Journal of Computational Physics
  357 (2018) 125--141.

\bibitem{regazzoni2019machine}
F.~Regazzoni, L.~Ded{\`e}, A.~Quarteroni, Machine learning for fast and
  reliable solution of time-dependent differential equations, Journal of
  Computational Physics 397 (2019) 108852.

\bibitem{regazzoni2020machine}
F.~Regazzoni, L.~Ded{\`e}, A.~Quarteroni, Machine learning of multiscale active
  force generation models for the efficient simulation of cardiac
  electromechanics, Computer Methods in Applied Mechanics and Engineering 370
  (2020) 113268.

\bibitem{hesthaven2018non}
J.~S. Hesthaven, S.~Ubbiali, Non-intrusive reduced order modeling of nonlinear
  problems using neural networks, Journal of Computational Physics 363 (2018)
  55--78.

\bibitem{ray2018artificial}
D.~Ray, J.~S. Hesthaven, An artificial neural network as a troubled-cell
  indicator, Journal of Computational Physics 367 (2018) 166--191.

\bibitem{VEMquality2021}
T.~Sorgente, S.~Biasotti, G.~Manzini, M.~Spagnuolo, The role of mesh quality
  and mesh quality indicators in the virtual element method, arXiv preprint
  arXiv:2102.04138.

\bibitem{dryden2016statistical}
I.~L. Dryden, K.~V. Mardia, Statistical shape analysis: with applications in R,
  Vol. 995, John Wiley \& Sons, 2016.

\bibitem{bishop2006pattern}
C.~M. Bishop, Pattern recognition and machine learning, springer, 2006.

\bibitem{ioffe2015batch}
S.~Ioffe, C.~Szegedy, Batch normalization: Accelerating deep network training
  by reducing internal covariate shift, in: International conference on machine
  learning, PMLR, 2015, pp. 448--456.

\bibitem{brenner2018virtual}
S.~C. Brenner, L.-Y. Sung, Virtual element methods on meshes with small edges
  or faces, Mathematical Models and Methods in Applied Sciences 28~(07) (2018)
  1291--1336.

\bibitem{beirao2017stability}
L.~Beir{\~a}o~da Veiga, C.~Lovadina, A.~Russo, Stability analysis for the
  virtual element method, Mathematical Models and Methods in Applied Sciences
  27~(13) (2017) 2557--2594.

\bibitem{droniou2021robust}
J.~Droniou, L.~Yemm, Robust hybrid high-order method on polytopal meshes with
  small faces, arXiv preprint arXiv:2102.06414.

\bibitem{mu2015shape}
L.~Mu, X.~Wang, Y.~Wang, Shape regularity conditions for polygonal/polyhedral
  meshes, exemplified in a discontinuous galerkin discretization, Numerical
  Methods for Partial Differential Equations 31~(1) (2015) 308--325.

\bibitem{petersen2020neural}
P.~C. Petersen, Neural network theory, University of Vienna.

\bibitem{kingma2014adam}
D.~P. Kingma, J.~Ba, Adam: A method for stochastic optimization, arXiv preprint
  arXiv:1412.6980.

\bibitem{talischi2012polymesher}
C.~Talischi, G.~H. Paulino, A.~Pereira, I.~F. Menezes, Polymesher: a
  general-purpose mesh generator for polygonal elements written in matlab,
  Structural and Multidisciplinary Optimization 45~(3) (2012) 309--328.

\bibitem{arnold2002unified}
D.~N. Arnold, F.~Brezzi, B.~Cockburn, L.~D. Marini, Unified analysis of
  discontinuous galerkin methods for elliptic problems, SIAM Journal on
  Numerical Analysis 39~(5) (2002) 1749--1779.

\bibitem{cockburn2012discontinuous}
B.~Cockburn, G.~E. Karniadakis, C.-W. Shu, Discontinuous Galerkin methods:
  theory, computation and applications, Vol.~11, Springer Science \& Business
  Media, 2012.

\bibitem{droniou2020interplay}
J.~Droniou, Interplay between diffusion anisotropy and mesh skewness in hybrid
  high-order schemes, in: International Conference on Finite Volumes for
  Complex Applications, Springer, 2020, pp. 3--23.

\end{thebibliography}

\end{document}